\documentclass[12pt]{article}
\usepackage{epsfig}
\usepackage{amsmath}
\usepackage{amsfonts,amsthm,amssymb}
\usepackage{xcolor}
\usepackage{float}
\usepackage{authblk}
\usepackage{indentfirst}
\usepackage[left=2cm,right=2cm,top=2cm,bottom=2cm]{geometry}
\usepackage{epsfig}
\usepackage{graphicx}
\usepackage{enumitem}
\usepackage{mathtools}

\newtheorem{lm}{Lemma}
\newtheorem{theorem}{Theorem}

\begin{document}

\begin{center}
{\Large \bf Universal bifurcation patterns in the unfolding of a pair of homoclinic tangencies}
\end{center}

\vspace{0.5cm}
\centerline{{\bf Sergey Gonchenko$^{1}$, Dongchen Li$^{2}$,  Dmitry Turaev$^2$}}

\vspace{0.5cm}
{$^1${Lobachevsky State University of Nizhny Novgorod,
Scientific and Educational Mathematical Center ``"Mathematics of Future Technologies'', Nizhny Novgorod, Russia} 

{$^2${Department of Mathematics, Imperial College London}}

\vspace{0.5cm}
\centerline{\footnotesize{\em e-mail: sergey.gonchenko@mail.ru; d.li@imperial.ac.uk; d.turaev@imperial.ac.uk;}}

\vspace{1.5cm}

\noindent{\bf Abstract.}
We study generic two- and three-parameter unfoldings of a pair of orbits of quadratic homoclinic tangency in strongly dissipative systems. We prove that the corresponding stability windows for periodic orbits have various universal forms: the so-called ``shrimps" (cross-road
areas), as well as spring and saddle areas, and (in three-parameter unfoldings)  ``pregnant shrimps'' -- specific types of transitions between the shrimps and spring or saddle areas. \\

\section{Introduction}

Numerous studies of chaotic dynamics in models of vastly different physical
nature reveal universal repetitive patterns in the bifurcation set in parameter space. These are
the so-called ``shrimps'' \cite{Gal93,Gal94,HGGYK99,BGU08,Lor08}, and also ``squids'' and ``cockroaches'' \cite{G85,DM00,GGO17}. This phenomenon was
discovered by C. Mira and coworkers in \cite{M87,BC91a,BC91b,BC91c} where less zoologically charged names were used -- ``cross-road area'',
``saddle area'', and ``spring area'', see Fig.~\ref{Donch_areas2}. It is well-known that chaos in
non-hyperbolic systems exhibits ``stability windows'' -- regions in parameter space which correspond
to emergence of stable periodic orbits. In the simplest case of strongly dissipative maps, the boundary of a stability window corresponds to a saddle-node bifurcation, where the stable periodic orbit is born, and
to a period-doubling bifurcation, where the periodic orbit becomes unstable.\footnote{We call a map strongly dissipative (or sectionally dissipative) if it contracts two-dimensional areas, so no orbit can have more than one positive Lyapunov
exponent. If this is not the case, then another stability boundary can exist which corresponds to the birth of an
invariant torus, leading to bifurcation patterns different from those discussed in this paper, see e.g. \cite{GOT12}. Note also
that in systems of differential equations (as opposed to maps) stability boundaries may also correspond to homoclinic loops or to a blue-sky catastrophe, even in the case of strong dissipation, see \cite{book}. We leave all these situations aside.}. The curious phenomenon which we discuss in this paper is that
quite often these two stability boundaries do not simply go parallel to each other (forming ''window streets'' by the terminology of \cite{Lor08}) but somehow conspire to form
non-trivial patterns mentioned above.

\begin{figure}[h]
\begin{center}
\includegraphics[width=16cm]{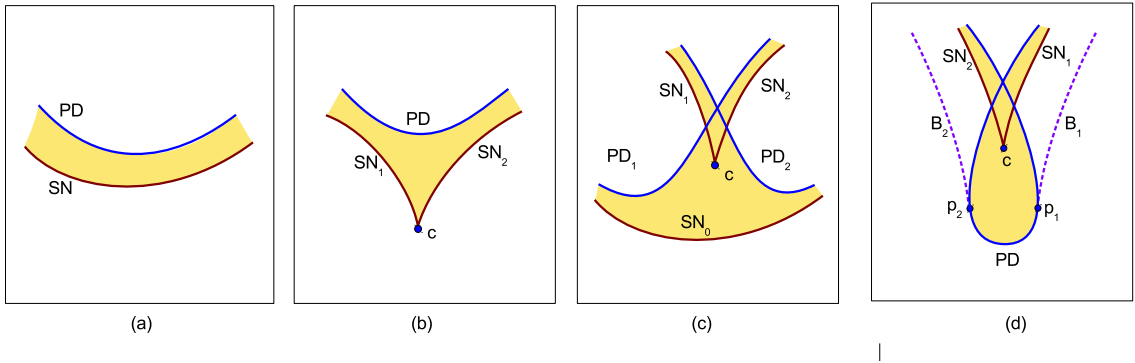}
\caption{{\footnotesize Examples of universal windows of stability: (a) ``band-like area'' called also ``stability window street''; (b) ``saddle area''  -- ``squid''; (c) ``cross-road area'' -- ``shrimp''; (d) ``spring area''  -- ``cockroach''. We denote by $SN$ and $PD$ curves of saddle-node (fold) and period-doubling (flip) bifurcations for fixed points and by $B$ a fold curve for period-2 points. The points  $c$ and $p$ correspond to codimension-2 bifurcations of fixed points -- pitch-fork (cusp point $c$)  and degenerate flip (with zero first Lyapunov value).}}
\label{Donch_areas2}
\end{center}
\end{figure}

We explain this by the ubiquity of homoclinic tangencies. Chaotic dynamics are in general associated with
hyperbolic sets and, in particular, to saddle periodic orbits and orbits of the intersection of their stable and unstable
manifolds, the homoclinics. If we observe changes in the structure of chaos when parameters of the system change,
it is natural to expect that the structure of the set of homoclinic orbits changes too, and non-transverse homoclinics appear at some parameter values, see Fig.~\ref{homtrtan}. According to the Newhouse theorem
\cite{N79,GST93b,PV94,R95}, this implies that  parameter values corresponding to
homoclinic tangencies are dense in some open regions in the parameter space. It is plausible that at least in the context of strongly dissipative maps these open regions
cover all parameter values corresponding to non-hyperbolic chaotic dynamics \cite{GST91a,GST93a,LGYK93,Palis2000}. It was discovered by Gavrilov and Shilnikov
\cite{GaS73} and Newhouse \cite{N74} that unfolding a homoclinic tangency in a strongly dissipative system is accompanied by the birth of stable periodic orbits. Thus, we can conclude that the particular stability windows that emerge due to bifurcations of homoclinic tangencies
should be universally present in chaotic models of arbitrary nature.

\begin{figure}[h]
\begin{center}
\includegraphics[width=16cm]{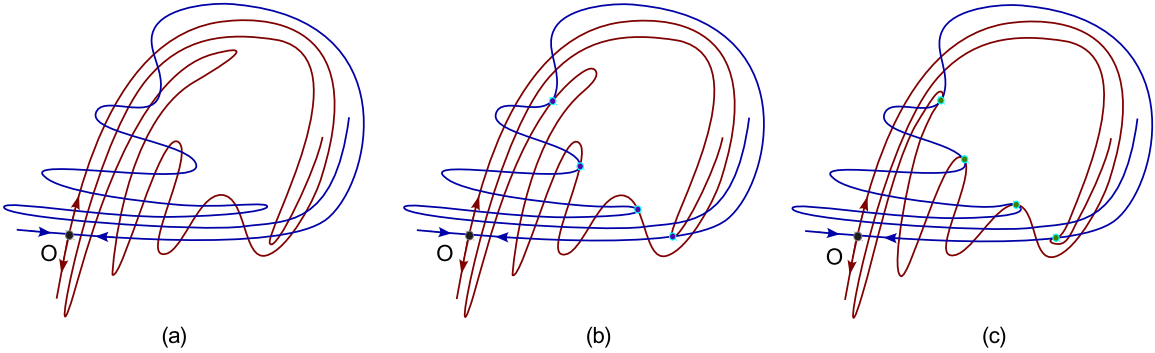}
\caption{{\footnotesize Homoclinic structures associated with a saddle fixed point $O$: (a) transverse homoclinics; (b) creation of a quadratic homoclinic tangency; (c) creation of a cubic homoclinic tangency (in general, two parameters are needed for this, see e.g. \cite{GST93a}).}}
\label{homtrtan}
\end{center}
\end{figure}

A generic one-parameter unfolding of a quadratic homoclinic tangency creates simple stability windows (the so-called ``window streets'' or ``band-like areas'') as in
Fig.~\ref{Donch_areas2}a . These are intervals of
parameter values bounded by a saddle-node bifurcation and a period-doubling bifurcation, see \cite{GaS73,G83}. However,
it was shown in \cite{GST91a,GST93a}, that bifurcations of a homoclinic tangency can never be completely captured by any
finite-parameter unfolding: increasing the number of independent parameters leads to new phenomena of a higher codimension, i.e., to phenomena that are not typically present in unfoldings governed by a smaller number of parameters. Thus, a two-parameter unfolding
of a quadratic homoclinic tangency creates cubic homoclinic tangencies \cite{GST93a,GST07}. The two-parameter stability windows near a cubic tangency were studied in \cite{G85} and they are exactly the ``saddle area'' (the ``squid'') and the ``spring area'' (the ``cockroach'') of Fig.~\ref{Donch_areas2}(b) and (d). It is also shown in \cite{GST93a} that a two-parameter unfolding of a single quadratic homoclinic tangency creates a pair of homoclinic tangencies. In the current paper we  study the stability windows due to the simultaneous breakup of two orbits of
homoclinic tangency and show that they have a form of a ``shrimp'' (the ``cross-road'') area.

We also study a generic three-parameter unfolding of a pair of orbits of quadratic homoclinic tangency and show that the corresponding stability windows
in the three-dimensional parameter space have a form of ``pregnant shrimps'' -- specific types of transitions between the cross-road, spring and saddle areas, as described in \cite{BC91a,BC91b,BC91c}. Thus, whenever we consider a model with chaotic behavior,
in the Newhouse regions in two-dimensional or three-dimensional parameter space there must exist, in abundance, stability windows of these particular shapes.

\section{Problem setting and results}\label{sec:setting,results}
\subsection{Simplest periodic orbits near a double homoclinic tangency}
Consider a $C^{r}$-smooth ($r\geq 2$) diffeomorphism $f$ of an $(n+1)$-dimensional $(n\geq 1)$ smooth manifold.
Let $f$ have a saddle periodic point $O$ with a one-dimensional unstable and an $n$-dimensional stable invariant manifolds $W^u(O)$ and $W^s(O)$. Assume that there exist two homoclinic orbits $\Gamma_1$ and $\Gamma_2$ such that $W^s(O)$ and $W^u(O)$ have {\em quadratic tangency} at the points of these two orbits (see Fig. \ref{neigh8}(a)).

\begin{figure}[!h]
\begin{center}
\includegraphics[width=14cm]{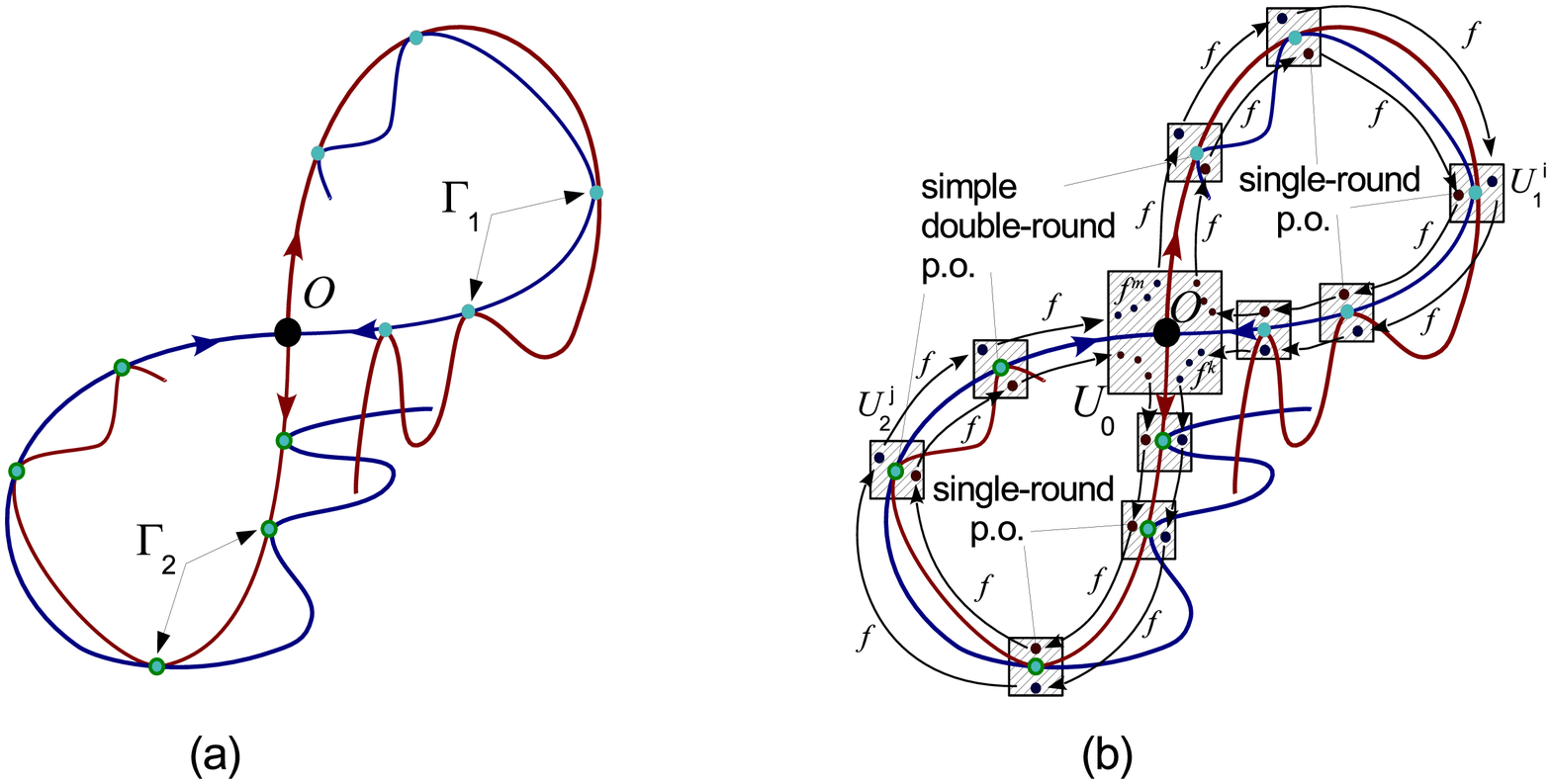}
\caption{{\footnotesize (a) A pair of homoclinic tangencies to a fixed point $O$ of a two-dimensional diffeomorphism;
(b) a neighborhood $U$ of $O\cup\Gamma_1\cup\Gamma_2$; here
$U = U_0\cup U_1\cup U_2$, where $U_0$ is a small disc containing $O$, the handles
$U_1 = \bigcup\limits_{k=1}^{n_1} U_1^k$ and $U_2 = \bigcup\limits_{k=1}^{n_2} U_2^k$ consist of small neighborhoods of those points of the orbits $\Gamma_1$ and
$\Gamma_2$ that do not belong to $U_0$. The iterations of $f$ for single-round and simple double-round periodic orbits are also shown schematically.}}
\label{neigh8}
\end{center}
\end{figure}

The closure of the homoclinic orbits $\Gamma_{1,2}$ is the union of $\Gamma_{1,2}$ and the orbit $L=\{f^j(O)\}_{j=0}^{\tau-1}$ of the point $O$
(where $\tau$ is a period of $O$, i.e., $f^\tau(O)= O$).
We take a small neighborhood $U$ of the set of these three orbits. Our goal is to study bifurcations of simplest periodic orbits lying in $U$.

Let $U_0$ be a small ball around $O$. Since the orbit $\Gamma_i$ ($i=1,2$) is homoclinic, only
finitely many points in $\Gamma_i$ lie outside of the neighborhood $\tilde U_0 = \bigcup_{j=0}^{\tau-1} f^j(U_0)$ of $L$. Namely, on $\Gamma_i$ there are points $P_i^+\in W^s_{loc}(O)\cap U_0$ and $P_i^-\in W^u_{loc}(O)\cap U_0$ such that $P_i^+=f^{N_i} P_i^-$ for some $N_i>0$ while $f^j(P_i^-) \not\in \tilde U_0$ for $j=1, \dots, N_i-1$.

Note that $U$ is a disjoint union of $\tilde U_0$ and a small neighborhoods $U_1$ and $U_2$ of $\Gamma_1\setminus \tilde U_0$ and,
respectively, $\Gamma_2\setminus \tilde U_0$ (the set $U_i$ is a disjoint union of small balls around the points $f^{j} P_i^-$, $j=1,\dots,p_i-1$,
see Fig. \ref{neigh8}(b)). As $L$ is hyperbolic, no other orbit can lie entirely in $\tilde U_0$.  So, any periodic orbit in $U$, other than $L$, must visit $U_1\cup U_2$. The same is true
for any map which is $C^1$-close to $f$: there is a uniquely defined hyperbolic continuation of $L$ in $\tilde U_0$ and any other orbit in $U$ must visit $U_1\cup U_2$.

A periodic orbit in $U$ is called {\em single-round} if it visits only $U_1$ or only $U_2$, and does it exactly once. A {\em simple double-round} periodic orbit visits $U_1$, gets into $\tilde U_0$, then leaves it and goes to $U_2$, and finally returns to $\tilde U_0$ and closes up.

Single-round orbits stay in a small neighborhood of only one of the orbits of homoclinic tangency, $\Gamma_1$ or $\Gamma_2$, so their bifurcations
are described by the theory of a homoclinic tangency, which starts with the work of Gavrilov and Shilnikov \cite{GaS73}, see \cite{G83,GST93c,GST08} for the multidimensional case. We, therefore, focus in this paper on the new case -- the simple double-round periodic orbits.

We call the map $T_0=f^\tau|_{U_0}$ the {\em local map} near the point $O$. We will also consider the {\em global maps} $T_1$ and $T_2$, where $T_i$ is the restriction of $f^{N_i}$ onto a sufficiently small ball around $P^-_i$, so it acts from a small neighborhood of $P^-_i$ to a small neighborhood of $P_i^+$.  The maps $T_0$ and $T_{1,2}$ remain well-defined for any $C^1$-small perturbation of $f$.

By definition, a single-round orbit corresponds to a fixed point of the first-return map $T_iT_0^k$ for some sufficiently large integer $k$ and $i=1$ or $i=2$.
The simple double-round orbits are defined as those corresponding to fixed points of the first-return map $T_{km}=T_2T_0^m T_1 T_0^k$ for some sufficiently large integers $m$ and $k$. We will call them $(k,m)$-{\em orbits}.

A very similar construction is applicable for dynamical systems with continuous time. Namely, we consider in this case a $C^{r}$-smooth ($r\geq 2$) system of differential equations or, equivalently, a  $C^{r}$-smooth vector filed $f$ on an $(n+2)$-dimensional $(n\geq 1)$ smooth manifold. Let $f$ have a saddle periodic orbit $L$, and let $O$ be a point of intersection of this periodic orbit with a small $(n+1)$-dimensional cross-section $U_0$.
Let $\tilde U_0$ be a sufficiently small neighborhood of $L$. Then, the {\em local  map} $T_0$ is defined as the map of first return to $U_0$ by the orbits of the flow defined by $f$, which lie in the neighborhood $\tilde U_0$.

We assume that the unstable manifold $W^u(O)$ is one-dimensional and the stable manifold $W^s(O)$ is $n$-dimensional and that $W^u(O)$ and $W^s(O)$ have two orbits of quadratic homoclinic tangency, $\Gamma_1$ and $\Gamma_2$. On $\Gamma_i$, we take points $P_i^+\in W^s_{loc}(O)\cap U_0$ and $P_i^-\in W^u_{loc}(O)\cap U_0$, and denote as $U_i$ a small neighborhood of the segment of $\Gamma_i$ between the points $P_i^-$ and $P_i^+$. The {\em global maps} $T_i$ ($i=1,2$) are defined as maps from a sufficiently small ball around $P^-_i$ in $U_0$ to a small ball around $P^-_i$ in $U_0$ along the orbits lying in $U_i$. These maps remain well-defined for any $C^1$-small perturbation of $f$.

A periodic orbit in $U=\tilde U_0\cup U_1\cup U_2$ is {\em single-round} if it visits only $U_1$ or only $U_2$, and does it exactly once.
A {\em simple double-round} periodic orbit visits $U_1$, gets into $\tilde U_0$, then leaves it and goes to $U_2$, and finally returns to $\tilde U_0$ and closes up.
By definition, a single-round orbit corresponds to a fixed point of the first-return map $T_iT_0^k$ in a small neighborhood of $P_i^+$ in $U_0$ for some sufficiently large integer $k$ and $i=1$ or $i=2$.
The simple double-round orbits $(k,m)$-{\em orbits} are defined as those corresponding to fixed points of the first-return map $T_{km}=T_2T_0^m T_1 T_0^k$ in a small neighborhood of $P_i^+$ in $U_0$ for some sufficiently large integers $m$ and $k$.

All the arguments in the proofs of the results below are done in terms of the local and global maps only and, with the above notations, are identical for the cases of discrete and continuous time (i.e., diffeomorphisms and smooth flows). We therefore adopt a neutral terminology in the formulation of the results, simply calling $f$ {\em a dynamical system}, implying that all the theorems hold true in both cases.

\subsection{Two-parameter unfoldings}
We start with analyzing bifurcations of the $(k,m)$-orbits in a generic two-parameter unfolding of the homoclinic tangencies of $f$. Consider a two-parameter family $f_{\varepsilon}$, which depends smoothly on $\varepsilon=(\varepsilon_1,\varepsilon_2)$, and let $f_0$ be the original system $f$. Introduce coordinates $(x,y)\in \mathbb{R}^n \times \mathbb{R}$ in $U_0$ such that the local stable and unstable manifolds are given for all small $\varepsilon$ by the equations $y=0$ and $x=0$, respectively. Define the {\em splitting parameter} $\mu_i$ for the tangency $\Gamma_i$ as the $y$-coordinate of the point $T_i P_i^-$ (see Lemma~\ref{lem:T_1,2}). The splitting parameters are
smooth functions of $\varepsilon$. The pair of tangencies is said {\em to unfold generically} if
\begin{equation}\label{eq:paraderi1}
\det \frac{\partial (\mu_1,\mu_2)}{\partial (\varepsilon_1,\varepsilon_2)} \neq 0.
\end{equation}
This will be our standing assumption; it allows us to take $\mu=(\mu_1,\mu_2)$ as new parameters, so we will use the notation $f_\mu$ from now on.

Denote by $\lambda_1,\dots,\lambda_n,\gamma$ the multipliers of $O$ (the eigenvalues of the derivative matrix of $T_0$ at $O$) ordered such that
$$|\gamma|  > 1 > |\lambda_1|  \ge\dots\ge |\lambda_n|.$$
Our main assumption on $f$ is the {\em strong dissipativity} (also called sectional dissipativity) condition
\begin{equation}\label{eq:strdis}
|\lambda_1\gamma|<1.
\end{equation}
In general, bifurcations of single-round orbits near a homoclinic tangency can be quite complicated \cite{GST08}. In the strongly dissipative case, however, the situation is simpler. It immediately follows from \cite{GaS73,G83} that in the $(\mu_1,\mu_2)$-parameter plane there are infinitely many bands
(``stability window streets'') $\sigma_k^1$ and $\sigma_k^2$ accumulating at the axes $\mu_1=0$ and $\mu_2=0$, respectively, as $k\to\infty$. The bands $\sigma_i^\alpha$ and $\sigma_j^\alpha$, where $\alpha=1,2$, do not intersect if $i\neq j$. The boundaries of $\sigma_k^\alpha$ and $\sigma_k^\alpha$ are curves corresponding to nondegenerate saddle-node ($SN$) and period-doubling ($PD$) bifurcations for orbits of period $k$, see Fig.~\ref{strips2}.

\begin{figure}[!h]
\begin{center}
\includegraphics[width=16cm]{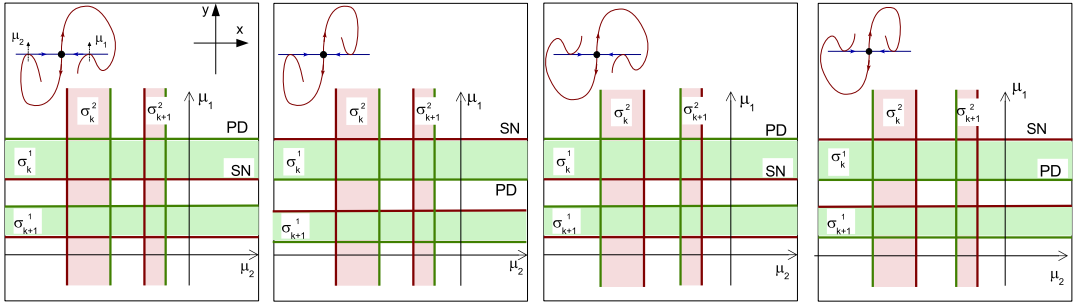}
\caption{{\footnotesize Bands $\sigma_k^1$ and $\sigma_k^2$ in the $(\mu_1,\mu_2)$-parameter plane for various types of double homoclinic tangencies in the case $\gamma>0$. Here we assume that the homoclinic $\Gamma_i$ split in such a way that $T_iW^u_{loc}(O)$ and $W^s_{loc}(O)$ do not intersect near $P_i^+$ for $\mu_i<0$ and intersect at two points at $\mu_i>0$.
Horizontal (green) and vertical (pink) bands correspond to stability window streets for single-round periodic orbits that lie, respectively, in $U_1$ and in $U_2$.}}
\label{strips2}
\end{center}
\end{figure}

The structure of stability windows for simple double-round orbits is different, as follows from Theorem \ref{thm:shr1} below.
Denote
\begin{equation}\label{thdf}
\theta= - \frac{\ln |\lambda_1|}{\ln|\gamma|}.
\end{equation}
Note that the strong dissipativity condition (\ref{eq:strdis}) implies $\theta>1$.

\begin{theorem}\label{thm:shr1}
Take any $\delta>0$. If $\theta>1$, and the two homoclinic tangencies are quadratic, then for every sufficiently large integers $m$ and $k$ satisfying
\begin{equation}\label{mkthet}
(\theta +\delta)^{-1} <\frac{m}{k} < \theta-\delta,
\end{equation}
there is a disc $\Delta_{k,m}$ in the $(\mu_1,\mu_2)$-plane and a ball $B_{k,m} \subset U_0$ such that the return map $T_{km}$ is defined
on $B_{k,m}$ for $\mu\in \Delta_{k,m}$ and, after an affine change of coordinates and parameters $R_{k,m}: (x,y,\mu_1,\mu_2)\mapsto (X,Y,M_1,M_2)$,
takes the form $T_{km}|_{B_{k,m}}: (X,Y) \mapsto (\bar X, \bar Y)$, where
\begin{equation}\label{eq:1mapsc}
\begin{array}{ll}
\bar X =o_{r}(1)_{k,m\to+\infty} ,\\
\bar Y = M_2 - (M_1 - Y^2)^2 +o_{r}(1)_{k,m\to+\infty}.
\end{array}
\end{equation}
Here
$o_{r}(1)_{k,m\to+\infty}$ stands for terms tending to zero along with all derivatives up to order $r$ as $k,m\to +\infty$. As $k,m\to+\infty$, the discs $\Delta_{km}$ converge to the origin $(\mu_1,\mu_2)=0$ and the balls $B_{km}$ converge to the homoclinic
point $P_2^+$ (see Fig.~\ref{frtms12n2}a). The image of $B_{km}\times \Delta_{km}$ by the rescaling map $R_{km}$ covers a centered at zero ball whose radius grows to infinity as the integers $k,m$ satisfying condition (\ref{mkthet}) tend to $+\infty$.
\end{theorem}

\begin{figure}[!h]
\begin{center}
\includegraphics[width=14cm]{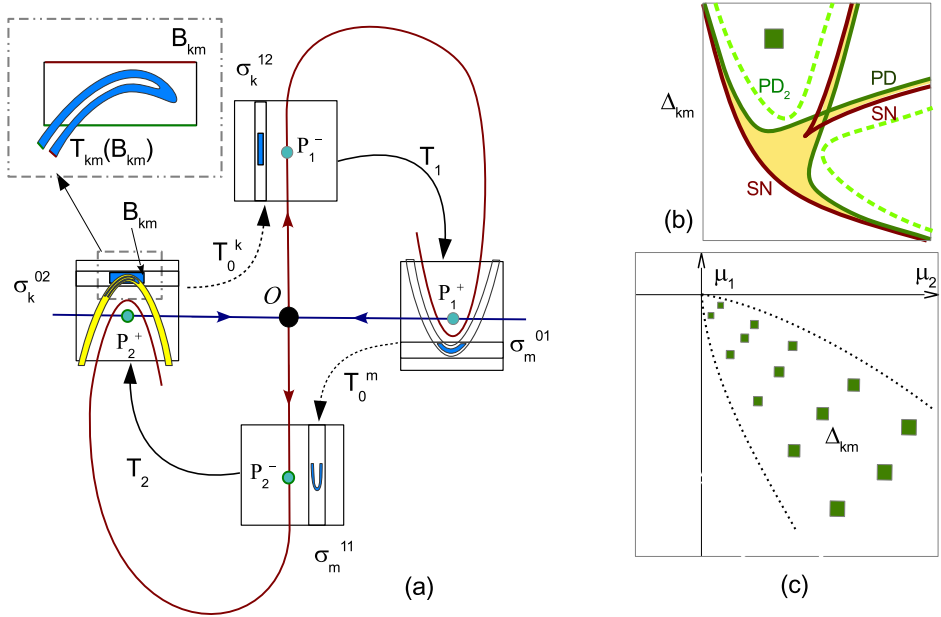}
\caption{{\footnotesize (a) The first-return map $T_{km}$ for simple double-round $(k,m)$-orbits; (b) the shrimp-like structure of a bifurcation set in $\Delta_{km}$;
(c) a cascade of the shrimp-like stability windows in the case of a double homoclinic configuration of Fig.(a).}}
\label{frtms12n2}
\end{center}
\end{figure}

This result (see the proof in Section \ref{sec:proofofthm}) implies that bifurcations of the
$(k,m)$-orbits for large $(k,m)$ satisfying condition (\ref{mkthet})
are the same as bifurcations of the fixed points of the one-dimensional map \footnote{Strictly speaking, it is true when $r\geq 5$ - because the map (\ref{2pa}) can undergo the codimension-2 flip (period-doubling) bifurcation where the first Lyapunov value vanishes and the second one is non-zero (the calculation of the second Lyapunov value requires at least five derivatives).}
\begin{equation}\label{2pa}
\bar Y = M_2 - (M_1 - Y^2)^2.
\end{equation}
Its bifurcations are well studied  \cite{M87,BC91a,BC91b,BC91c}; the bifurcation curves corresponding to the saddle-node and period-doubling bifurcations of the
fixed points form the ``shrimp'' (the ``cross-road area'') pattern as in Fig.~~\ref{frtms12n2}b. Thus, Theorem~\ref{thm:shr1} shows that
\begin{itemize}
\item
{\em in generic two-parameter families of strongly dissipative systems, points corresponding to the existence of a pair of homoclinic tangencies are accumulated by a cascade of the shrimp-like stability windows, see Fig.~\ref{frtms12n2}c.}
\end{itemize}
Of course, one does not need to be restricted by bifurcations of fixed points only: by Theorem~\ref{thm:shr1}, whichever non-degenerate bifurcations of codimension 1 and 2 are present in map (\ref{2pa}), they are repeated infinitely many times in any generic two-parameter unfolding of
a pair of homoclinic tangencies of a strongly dissipative system.\\

The idea of rescaling the first-return map was first proposed in \cite{LTY} for the study of a single quadratic homoclinic tangency.
For a single $n$-th order homoclinic tangency, the first-return map  in a generic $(n-1)$-parameter unfolding
gets, after a rescaling, arbitrarily close to the polynomial map
$$\bar Y = \sum\limits_{i=0}^{n-2} M_{i+1} Y^i \pm Y^n,$$
where we can always take the minus sign  in front of $Y^n$ when $n$ is even, see \cite{GST96,G96}. For $n=2$ (quadratic tangency), this gives us the  the parabola map
$$\bar Y = M_1 - Y^2;$$
a cubic homoclinic tangency gives rise to the maps
$$\bar Y = M_1 + M_2 Y \pm Y^3.$$

\begin{figure}[!h]
\begin{center}
\includegraphics[width=14cm]
{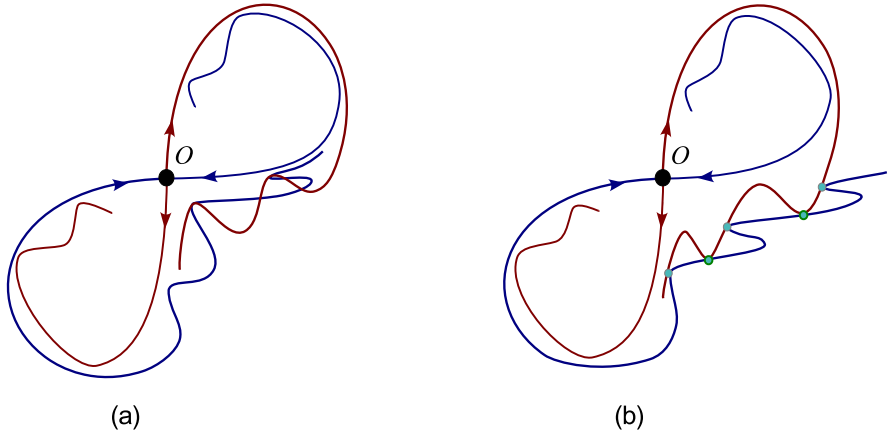}
\caption{{\footnotesize A schematic picture at perturbations of a diffeomorphism from Fig.~\ref{neigh8}(a) of  (a) a double-round cubic homoclinic tangency and (b) a pair of double-round  homoclinic tangencies. }}
\label{cubdouble}
\end{center}
\end{figure}

The stability windows for the parabola map and the cubic maps are shown in Figure~\ref{Donch_areas2}: the ``window streets'' appear in the parabola map; the so-called ``saddle area'' (``squid'') and the ``spring area'' (``cockroach'') patterns appear in the cubic maps with, respectively, the ``+'' and ``-'' sign \cite{G85,GST96}. Along with the shrimp pattern produced by the double parabola map (\ref{2pa}), they form the 4 universal bifurcation patterns
that appear in generic two-parameter families of strongly dissipative systems from the Newhouse domain.
In particular, one can show (cf. \cite{GST93a,GST99}) that the splitting of two quadratic homoclinic tangencies creates points in the
$(\mu_1,\mu_2)$-plane that correspond to a cubic tangency, see Fig.~\ref{cubdouble}b. So all the 4 patterns are present in a generic two-parameter family $f_\mu$ which we consider here.
For example, we see these patterns in the parameter plane of a periodically perturbed autonomous system with a homoclinic figure-eight.
Thus, in \cite{GSV13} the existence of infinite cascades of parameter values corresponding to cubic and pairs of quadratic homoclinic tangencies
was proven, and cascades of stability windows of of different types were found, see Fig.~\ref{fig8hom}.

\begin{figure}[!h]
\begin{center}
\includegraphics[width=14cm]
{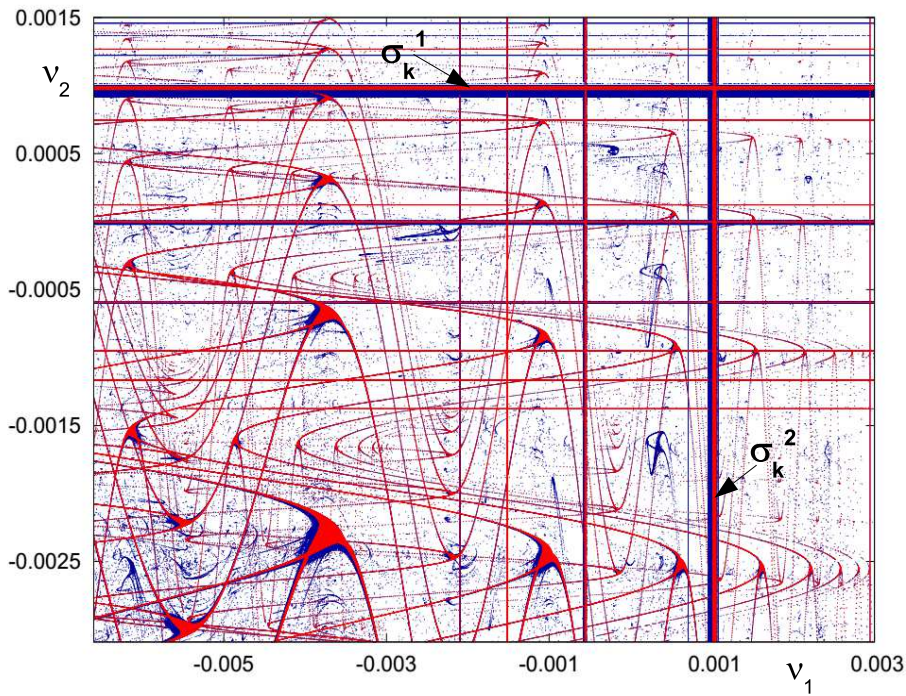}
\caption{{\footnotesize Bifurcation sets in the $(\nu_1,\nu_2)$-parameter plane
for a  periodically perturbed system with a homoclinic figure-eight; red stability windows correspond to single-round periodic orbits and blue windows to multi-round periodic orbits). }}
\label{fig8hom}
\end{center}
\end{figure}

\subsection{Three-parameter unfoldings and the main theorem}
By itself, the result of Theorem \ref{thm:shr1} is not surprising: as a single round near a homoclinic tangency scales to the parabola map, it is natural
that passing along two homoclinic tangencies gives rise to the composition of the parabola maps. A similar result can be extracted e.g. from \cite{Be18}. However, we want to stress that bifurcations of $(k,m)$-orbits are, in fact, richer than what happens in map (\ref{2pa}). Our main result, Theorem \ref{thm:shr1-3}, shows that the bifurcation patterns for the $(k,m)$-orbits which violate condition (\ref{mkthet}) are generated by the three-parameter
family of maps
\begin{equation}\label{eq:intro:shrimp2}
\bar Y = M_2 - (M_1-Y^2)^2 + M_3Y.
\end{equation}

Accordingly, the description of the stability windows for such orbits requires the analysis of a {\em three-parameter unfolding} of the map $f$. The choice of the third control parameter depends on whether the point $O$ is a saddle or a saddle-focus.
We assume that $\lambda_1$, the largest (in the absolute value) of the multipliers inside the unit circle, is simple. It is either real or complex; in the latter case $\lambda_1=\lambda_2^*$. We assume that
$$|\gamma|  > 1 > |\lambda_1|  > |\lambda_j| \qquad \mbox{for} \;\; j>1$$
if $\lambda_1$ is real --  then $O$ is called {\em a saddle}, or
$$|\gamma|  > 1 > \lambda  > |\lambda_j| \qquad \mbox{for} \;\; j>2$$
if $\lambda_{1,2}=\lambda e^{\pm i\varphi}$ where $\varphi\in (0,\pi)$  --  then $O$ is called {\em a saddle-focus}.

We consider $C^r$-families of systems (diffeomorphisms or flows) $f_{\varepsilon}$, $\varepsilon=(\varepsilon_1,\varepsilon_2,\varepsilon_3)$,
where $f_0$ is the original system $f$ with a pair of orbits of homoclinic tangency. The splitting parameters $\mu_{1,2}$, the ratio $\theta$ defined by (\ref{thdf}) and the argument $\varphi$ of $\lambda_1$ (when $\lambda_1$ is complex) are smooth functions of $\varepsilon$.
We put $\mu_3 = \theta(\varepsilon) - \theta(0)$ if $O$ is a saddle, and $\mu_3 = \varphi(\varepsilon) - \varphi(0)$ if $O$ is a saddle-focus. Assume that
\begin{equation}\label{eq:paraderi2}
\det\frac{\partial (\mu_1,\mu_2,\mu_3)}{\partial (\varepsilon_1,\varepsilon_2,\varepsilon_3)} \neq 0,
\end{equation}
which allows us to write the family as $f_{\mu}$ with $\mu=({\mu_1,\mu_2,\mu_3})$.

Recall that $\theta$ and $\varphi$ are continuous invariants (moduli) of topological equivalence of systems with homoclinic tangencies: if two systems have
different values of $\theta$ or $\varphi$, they are not topologically equivalent \cite{GS87,G89,GS92}. In some cases they are also moduli of the $\Omega$-equivalence (topological equivalence on the set of nonwandering orbits) \cite{GS90,GST96,G00,G02}. This means that any change in the values of $\theta$ or $\varphi$ leads to bifurcations of non-wandering orbits, e.g. periodic orbits, even when the tangency is not split \cite{GS86,GS87,GST99}. Therefore, using the moduli
as control parameters is natural in the study of homoclinic bifurcations \cite{GST91b,GST93a,GST93c,GStT96,GST96,GSStT97,GST07,GST08}.

The last assumption we need for Theorem \ref{thm:shr1-3}, is that both homoclinic tangencies $\Gamma_1$ and $\Gamma_2$ are {\em simple}. The simplicity condition is a version of the quasi-transversality condition from \cite{NPT}. For two-dimensional maps with a saddle and three-dimensional maps with a saddle-focus, every quadratic homoclinic tangency is automatically simple. In higher dimensions, the stable manifold $W^s_{loc}(O)$
contains the  strong stable invariant submanifold $W^{ss}_{loc}(O)$. This submanifold is  $C^r$-smooth and  is uniquely defined by the condition that
it is an invariant manifold tangent at $O$ to the eigenspace (of the linearization of $T_0$ at $O$) corresponding to the multipliers
$\lambda_2,\dots, \lambda_n$ in the case of saddle, and $\lambda_3,\dots, \lambda_n$ in the case of saddle-focus. We denote $n_s=1$ if $O$
is a saddle, and $n_s=2$ if $O$ is a saddle-focus. Then the strong-stable manifold is $(n-n_s)$-dimensional.
It is also well-known (see e.g. \cite{book}) that $W^s_{loc}$ carries a uniquely
defined $C^{r}$-smooth strong-stable invariant foliation $\mathcal{F}^{ss}$
consisting of $(n-n_s)$-dimensional leaves; and $W^{ss}_{loc}$ is the leaf of $\mathcal{F}^{ss}$ which
contains $O$. Another fact we use (see, for example, \cite{HPS,NPT,T96,book})
is that the invariant unstable manifold $W^{u}(O)$ is a part
of the extended unstable manifold $W^{ue}$ which is an $(n_s+1)$-dimensional invariant manifold, tangent at
$O$ to the eigenspace corresponding to the eigenvalues $\gamma$ and $\lambda_1$ if $O$ is a saddle, or
$\gamma$ and $\lambda_1$, $\lambda_2$, if $O$ is a saddle-focus. The manifold $W^{ue}$ is at least
$C^{1+\epsilon}$. It is not uniquely defined, but it always contains $W^u_{loc}$ and any two such manifolds are
tangent to each other at the points of $W^u_{loc}$ (see Fig. \ref{fig:nonsimple}).

We can now interpret the simplicity condition as

\noindent\textbf{(C1)} The homoclinic orbits $\Gamma_{i}$ $(i=1,2)$ do not have points in the strong-stable manifold (i.e., $M_{i}^+\not\in W^{ss}_{loc}(O)$).\\
\noindent\textbf{(C2)} The extended unstable manifold is transverse to the strong-stable foliation at the points of $\Gamma_1$ and $\Gamma_2$, that is, $T_i W^{ue}_{loc}$ is transverse at the point $P_i^+$ to the leaf of $F^{ss}$ which passes though this point. (See Fig.~\ref{fig:nonsimple} for an illustration).

\begin{figure}[h!]
\begin{center}
\includegraphics[width=0.85\columnwidth]
{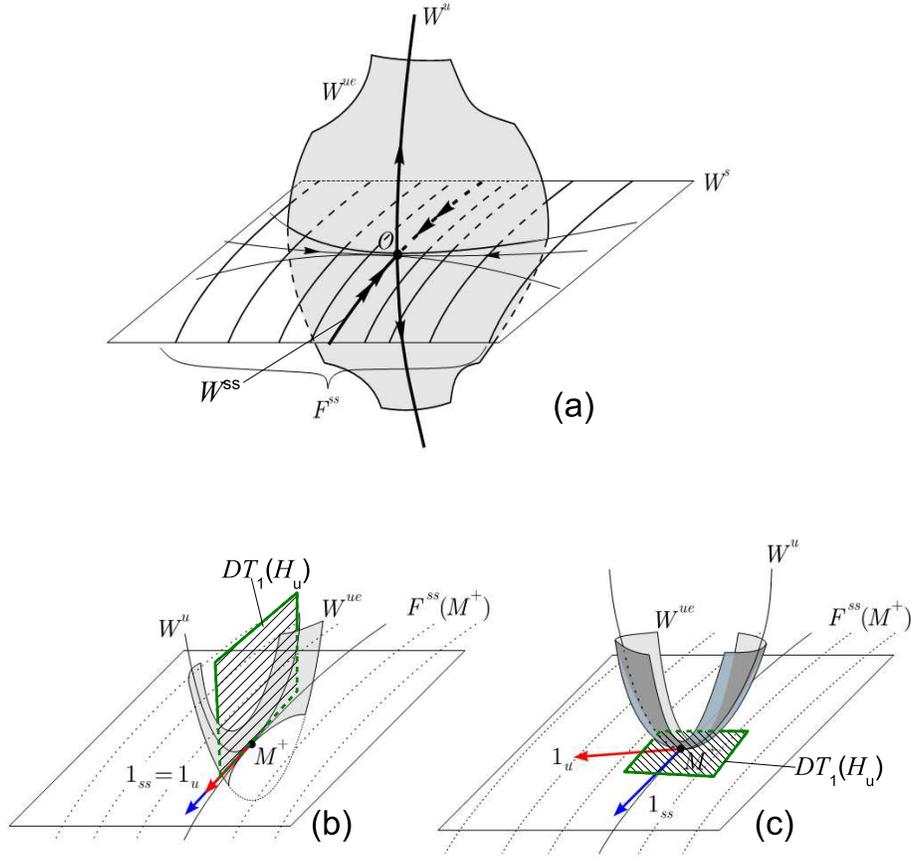} \caption{(a) A geometric illustration of the invariant manifolds ($W^u,W^s,W^{ss},W^{ue}$) and foliation $\mathcal{F}^{ss}$ existing near the fixed point $O$; (b) and (c) are examples of two types of non-simple homoclinic tangencies which correspond two different cases of nontransversality of $W^{ue}$ and $\mathcal{F}^{ss}(P^+)$. }
\label{fig:nonsimple}
\end{center}
\end{figure}

\begin{theorem}\label{thm:shr1-3}
If the system is strongly dissipative and the two tangencies are quadratic and simple (i.e., C1 and C2 are satisfied), then there exists a sequence of integers $k,m \to +\infty$ such that for every sufficiently large $k$ and $m$ from this sequence,
there exists a region $\Delta_{km}$ in the $(\mu_1,\mu_2,\mu_3)$-space and a ball $B_{km} \subset U_0$ such that
after a change of coordinates and parameters
$R_{km}: (x,y,\mu_1,\mu_2,\mu_3)\mapsto (X,Y,M_1,M_2,M_3)$, which is a composition of a transformation independent of $k$ and $m$
and an affine transformation with $k,m$-dependent coefficients, the return map $T_{km}$
on $B_{km}$ for $\mu\in \Delta_{km}$
takes the form $T_{km}|_{B_{km}}: (X,Y) \mapsto (\bar X, \bar Y):$
\begin{equation}\label{eq:2mapsc}
\begin{array}{ll}
\bar X_1 = M_1 - Y^2 + o_{r}(1)_{k,m\to+\infty},\\
\bar X_2 = o_{r}(1)_{k,m\to+\infty},\\
\bar Y = M_2 - (M_1 - Y^2)^2 +M_3 Y +o_{r}(1)_{k,m\to+\infty}.
\end{array}
\end{equation}
As $k,m\to+\infty$, the balls $\Delta_{km}$ converge to the origin $\mu=0$ and the balls $B_{km}$ converge to the homoclinic
point $P_1^+$. If $O$ is a saddle-focus, the image of $B_{km}\times \Delta_{km}$ by the rescaling $R_{km}$ covers a centered at zero ball whose radius grows to infinity as $k,m\to+\infty$. If $O$ is a saddle, then the image of $B_{km}\times \Delta_{km}$ by $R_{km}$ covers the intersection of a centered at zero ball whose radius grows to infinity as $k,m\to+\infty$ and either the set $M_3> \delta_{km}$
or the set $M_3 < - \delta_{km}$ where $\delta_{km}$ is a converging to zero sequence of positive numbers.
\end{theorem}

The proof is given in Section~\ref{sec:proofofthm}, along with the formulas for the transformation $R_{km}$ that relates the original
coordinates and parameters and the rescaled ones. Note that non-zero finite values of $M_3$ correspond, in the case of saddle, to
$\lambda_1^k\gamma^m$ (if $m>k$) or $\lambda_1^m\gamma^k$ (if $m<k$) being bounded away from zero and infinity. This implies
$m/k \to \theta$ or, respectively, $k/m \to \theta$ as $k,m\to +\infty$, so Theorem
\ref{thm:shr1-3} gives, in the saddle case, the first-return map for orbits corresponding to the $(k,m)$ values at the border of
those described by Theorem \ref{thm:shr1} (see condition (\ref{mkthet})). In the saddle-focus case there are no restrictions on $m/k$, as the parameter $M_3$ strongly depends on $\varphi$. Note also that in the saddle case $M_3$ has a definite sign for a given $(k,m)$; it can be different for another pair $(k,m)$, depending on the sign of $\lambda_1$ and $\gamma$ and the coefficients of the global maps $T_{1,2}$, see the proof of Theorem \ref{thm:shr1-3} for details.

\begin{figure}[!htb]
\begin{center}
\includegraphics[width=16cm]{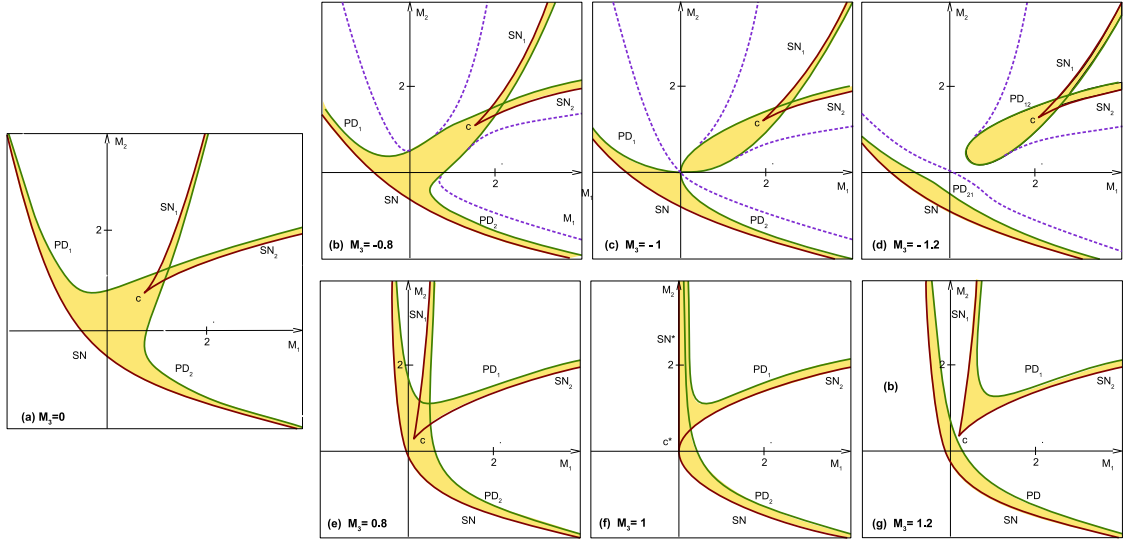}
\vspace{-0.5cm}
\caption{{\footnotesize An example of transition, in map (\ref{eq:intro:shrimp2}), from a ``shrimp'' to a ``cockroach'' when $M_3<0$, based on \cite{BC91a}; and from ``shrimp'' to ``squid'' when $M_3>0$, based on \cite{BC91c}. Here we show bifurcation curves for fixed points, saddle-node (SN) and period-doubling (PD) ones, as well as the saddle-node bifurcation curves (dashed lines) for period-2 points. (a)``Shrimp'' in the $(M_1,M_2)$-plane; (b) ``pregnant shrimp''; (c) the point $(M_1=0,M_2=0)$ corresponds to a codimension-3 bifurcation, when the fixed point $Y=0$ has the multiplier $-1$ and zero first and second Lyapunov values (the map becomes $\bar Y = Y - Y^4$); (d)``cockroach'';  (e)``shrimp''; (f) the point $(M_1=0,M_2=0)$ corresponds to a codimension-3 bifurcation, when the fixed point $Y=0$ has the multiplier $1$ and zero first and second Lyapunov values (the map becomes $\bar Y = - Y - Y^4$); (g) ``squid'': the curves $SN$ and $PD_2$ are separated.}}
\label{Donch2}
\end{center}
\end{figure}

By Theorem \ref{thm:shr1-3}, the generic three-parameter unfolding of the pair of simple quadratic homoclinic tangencies
exhibits a cascade of stability windows whose shape is the same as in the three-parameter family (\ref{eq:intro:shrimp2}) (possibly restricted
to only positive or negative values of $M_3$). This family was studied e.g in \cite{BC91a,BC91c}. The change in the structure of the stability windows in the $(M_1,M_2)$-parameter plane when $M_3$ varies is shown in Fig.~\ref{Donch2}. Note that the cases $M_3<0$ and $M_3>0$ are different. When $M_3<0$, one can see the transition from a ``shrimp'' to a ``pregnant shrimp'' to a ``cockroach''; when $M_3>0$, the ``shrimp'' - ``squid'' transition happens. A detailed analysis of these and other transitions was done in \cite{BC91a,BC91b,BC91c,OMCMML19}.

\section{Normal forms of the local and global maps}\label{sec:normalforms}
Let $f_{\varepsilon}$ be a generic unfolding family of $C^r$ diffeomorphisms ($r\geq 2$), where $f_0 = f$ is strongly dissipative and has a pair of orbits of quadratic homoclinic tangency. As mentioned before, we will consider for Theorem \ref{thm:shr1} the two-parameter unfolding with $\mu=(\mu_1,\mu_2)$ consisting of two splitting parameters, and for Theorem \ref{thm:shr1-3} the three-parameter unfolding with an additional parameter $\mu_3 = \theta(\varepsilon) - \theta(0)$ if $O$ is a saddle, and $\mu_3 = \varphi(\varepsilon) - \varphi(0)$ if $O$ is a saddle-focus. Due to conditions \eqref{eq:paraderi1} and \eqref{eq:paraderi2}, we can just denote the family by $f_{\mu}$.

Recall that the local map is defined as $T_0\equiv f^{\tau}:U_0\to U_0$, where $U_0$ is a small neighborhood of $O$; we take points $P_i^+\in W^s_{loc}(O)\cap \Gamma_i$ and $P_i^-\in W^u_{loc}(O)\cap \Gamma_i$ inside $U_0$ such that $P_i^+=f^{N_i} P_i^-$ for some $N_i>0$ while $f^j(P_i^-) \not\in \tilde U_0$ for $j=1, \dots, N_i-1$. Now fix in $U_0$ some small neighborhoods $\Pi_i^+\ni P_i^+$ and $\Pi_i^-\ni P_i^-$ satisfying $T_0(\Pi_i^{+})\cap\Pi_i^{+}=\emptyset$ and
$T_0^{-1}(\Pi_i^{-})\cap\Pi_i^{-}=\emptyset$. The global maps are $T_{1}\equiv f^{N_1} :\Pi_1^-\to\Pi_1^+$ and $T_{2}\equiv
f^{N_2} :\Pi_2^-\to\Pi_2^+$, which are diffeomorphisms into their ranges.

Observe that the iterates $T^k_0(\Pi_i^{+})$ intersect $\Pi_j^{-}$ ($i,j=1,2$) for every sufficiently large $k$. The domain of $T^k_0: \Pi_i^{+}\to \Pi_j^-$ consist of infinitely many pairwise disjoint strips $\sigma_{k}^{0ij}\subset \Pi_i^+$ which accumulate on $W^s_{\mathrm{loc}}(O)\cap\Pi_i^+$ as $k\to+\infty$. Its range is formed by infinitely many strips $\sigma_{k}^{1ij}=T_{0}^k(\sigma_k^{0ij})\subset\Pi_i^-$ accumulating on $W^u_{\mathrm{loc}}(O)\cap\Pi_i^-$ as $k\to+\infty$ (see Fig. \ref{fig:locmaps}).

\begin{figure}[h!]
\begin{center}
\includegraphics[width=16cm]%, height=7cm]
{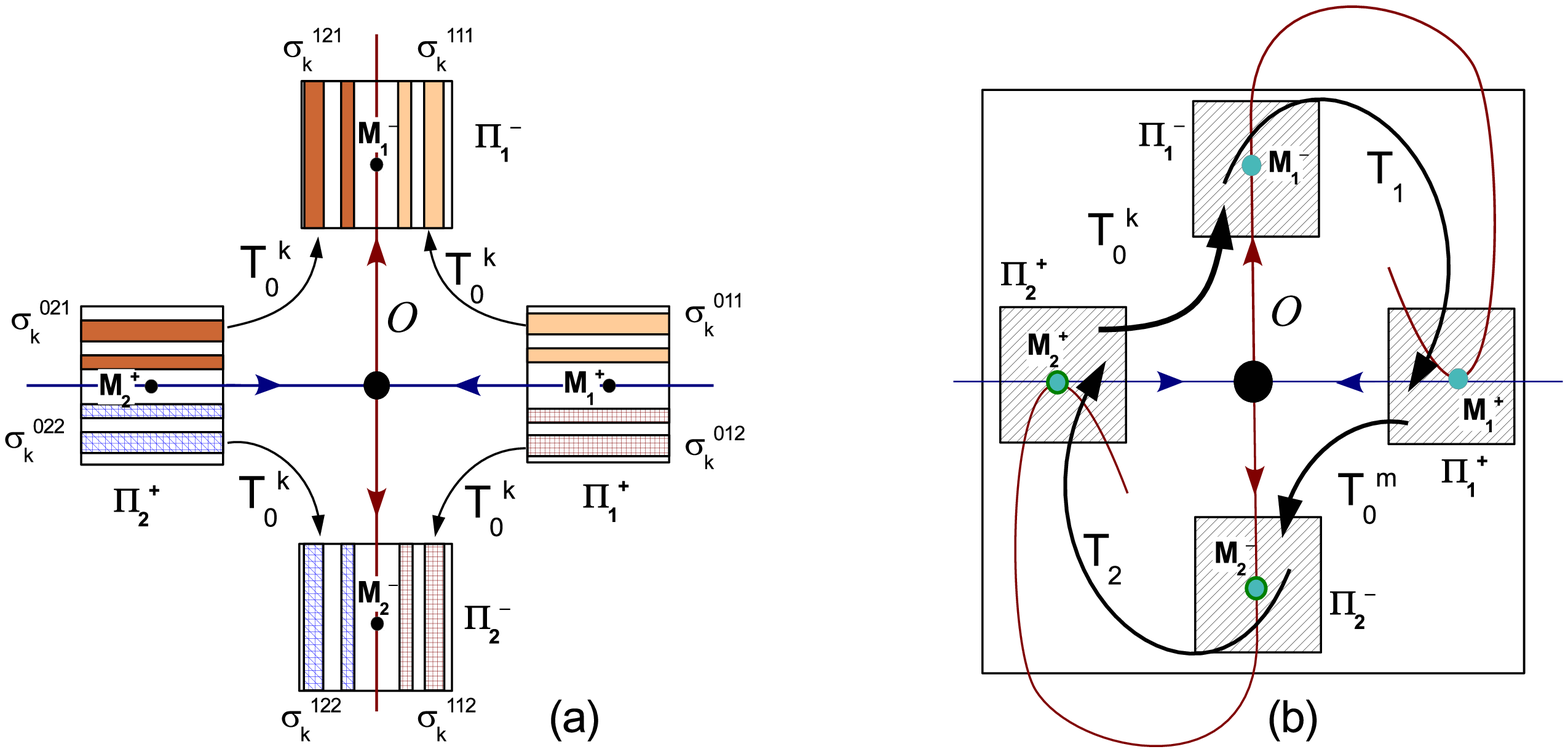} \caption{ (a) The action of iterations of $T_{0}$ from $\sigma^{0ij}\in\Pi^+_i$ to $\sigma^{1ij}\in\Pi^-_j$. (b) A geometric structure of the homoclinic points $M_1^+,M_1^-,M_2^+$ and $M_2^-$ and their neighborhoods. Schematic actions of the first return maps $T_{km}= T_2\circ T_0^m\circ T_1\circ T_0^k$.
}
\label{fig:locmaps}
\end{center}
\end{figure}

For the family $f_\varepsilon$,
one can always consider such $C^r$-coordinates in $U_0$
that the fixed point $O_\varepsilon$ of the local map $T_0$ is at
the origin for all small $\varepsilon$, and hence we drop the subscription of $O$. Moreover, we can write
$T_0(\varepsilon)$ in the form
\begin{equation}\label{eq:normalform:1}
\begin{array}{l}
\bar{x}=A_1x+\dots,\; \bar{u}=A_2u+\dots, \; \bar{y}=\gamma y+\dots,
\end{array}
\end{equation}
where $x\in \mathbb{R}^{n_s},y\in \mathbb{R}^{1},u\in \mathbb{R}^{m-n_s}$, the
dots stand for nonlinear terms, and the eigenvalues of the matrices
$A_1$ and $A_2$ are the stable leading
multipliers and the stable non-leading
multipliers, respectively. We say that $x$ are
the leading (stable) coordinates and  $u$ are
the non-leading ones. Note that if $\lambda_1$ is real, then
$A_1=\lambda_1$ and $x$ is a scalar; and if $\lambda_1$ is complex, then
$\lambda_1 = \bar\lambda_2 = \lambda e^{\pm i\varphi}$, $x=(x_1,x_2)$
and
$A_1=\lambda\left(
\begin{array}{cc}
\cos \varphi & \; -\sin \varphi\\
\sin \varphi & \; \cos \varphi
\end{array}
\right)$.

It is well-known that there is a $C^r$ coordinate transformation straightening the local stable and unstable manifolds of $O$
so that they acquire equations $W^s_{loc}=\{y=0\}$ and $W^u_{loc}=\{x=0,u=0\}$ and the
map $T_0$ takes, locally, the form
\begin{equation}
\begin{array}{l}
\bar x = A_{1}(\varepsilon)x + g_{1}(x,u,y,\varepsilon),\qquad
\bar u = A_{2}(\varepsilon)u + g_{2}(x,u,y,\varepsilon), \\
\bar y = \gamma(\varepsilon)y + h(x,u,y,\varepsilon),
\end{array}
\label{str}
\end{equation}
where the nonlinearities $g_{1,2}$ and $h$
vanish at the origin along with their first derivatives, and, moreover,
\begin{equation}
g(0,0,y,\varepsilon) \equiv  0, \qquad h(x,u,0,\varepsilon) \equiv  0
\label{pqsu}
\end{equation}
for $(x,u,y,\varepsilon)$ small.

Bringing the local map to the form (\ref{str}) is not enough for our purposes
because, for a general choice of the functions $g$ and $h$, the iterations
$T_0^k$ can deviate too much from its linear part. Essentially, this
means that the right-hand side of (\ref{str}) contains ``too many''
non-resonant terms. Fortunately, infinitely many of nonresonant terms can be
eliminated by means of some additional smooth transformation of coordinates,
as the following lemma from \cite{GST08} shows. This lemma is a generalization of similar results
of \cite{Afr84,GS90,GS92,book} and allows us to achieve a higher (i.e., $C^r$) smoothness of the coordinate transformation.

\begin{lm} \cite{GST08}
At all sufficiently small $\varepsilon$, there
exists a local $C^{r}$-transformation of coordinates after which
the map $T_0(\varepsilon)$ keeps its form {\rm (\ref{str}),(\ref{pqsu})} while
the functions $p$ and $q$ satisfy additional identities
\begin{equation}
g_{1}(x,u,0,\varepsilon) \equiv  0, \qquad
h(0,0,y,\varepsilon) \equiv  0,
\label{pq1}
\end{equation}
\begin{equation}
\frac{\partial g}{\partial x}(0,0,y,\varepsilon) \equiv  0, \qquad
\frac{\partial h}{\partial y}(x,u,0,\varepsilon) \equiv  0.
\label{pq12}
\end{equation}
\label{lem1}
\end{lm}

When the map $T_0(\varepsilon)$ is written in the coordinates of Lemma~\ref{lem1},
we say that it is in {\it the main normal form}. The identities (\ref{pq1}) imply that map $T_0$ becomes both linear on $W^u_{loc}=\{x=0,u=0\}$, i.e. the restriction $T_0|_{W^u_{loc}}$ has a form $\bar y = \gamma(\varepsilon)y$, and linear in $x$ on $W^s_{loc}=\{y=0\}$, i.e. the restriction $T_0|_{W^s_{loc}}$ has a form
$\bar x = A_{1}(\varepsilon)x,\; \bar u = A_{2}(\varepsilon)u + g_{2}(x,u,0,\varepsilon)$.
The latter also means that the strong stable foliation $\mathcal{F}^{ss}$ is straightened, i.e., its leaves take the form $\mathcal{F}^{s}:\{x=const,y=0\}$.
It is important that when $T_0$ is brought to this
normal form, the iterations $T_0^k\;:\;U_0\to U_0$ of the local map
do not differ too much from the iterations of the linearized map
at all large $k$. Let
$(x_k,u_k,y_k)=T_{0}^k(x_0,u_0,y_0)$. It has been known since
\cite{Sh67} (see also \cite{Sh68,AS73}) that $(x_k,u_k,y_0)$ are
uniquely defined functions of
$(x_0,u_0,y_k)$ for any $k\geq 0$.

\begin{lm} \cite{GST08}
When the local map $T_0$ is brought to the main normal form,
the following relations hold for all small $\varepsilon$ and all large $k$:
\begin{equation}
\begin{array}{l}
x_k \;-\; A_1^{k}(\varepsilon) x_0 \;=\;
\hat\lambda^{k}g_{k}(x_0,u_0,y_k,\varepsilon), \qquad
u_k \;=\; \hat\lambda^{k}\hat g_{k}(x_0,u_0,y_k,\varepsilon),\\
y_0 \;-\; \gamma^{-k}(\varepsilon) y_k \;=\;
\hat\gamma^{-k}h_k(x_0,u_0,y_k,\varepsilon),
\end{array}
\label{T0kk}
\end{equation}
where $\hat\lambda$  and $\hat\gamma$ are some constants such that
$0 < \hat\lambda < \lambda \;,\; \hat\gamma > \gamma$, and the functions
$g_{k},h_k,\hat g_k$ are uniformly bounded for all
$k$, along with the derivatives up to the order $(r-2)$.
The derivatives of order $(r-1)$ are estimated as
\begin{equation}\label{cr-1}
\left\| \frac{\partial^{r-1} (x_k - A_1(\varepsilon)^k x_0, \; u_k)}
{\partial^{r-1} (x_0,u_0,y_k,\varepsilon)}\right\|
=o(\|\lambda_1\|^k), \qquad
\left\| \frac{\partial^{r-1} (y_0 - \gamma(\varepsilon)^{-k} y_k)}
{\partial^{r-1} (x_0,u_0,y_k,\varepsilon)}\right\|
=o(\|\gamma\|^{-k}),
\end{equation}
while the derivatives of order $r$ are estimated as
\begin{equation}\label{cre}
\|x_k,u_k,y_0\|_{_{C^{r}}}=o(1)_{k\rightarrow\infty};
\end{equation}
these estimates do not include derivatives with more than $(r-2)$
differentiations with respect to $\varepsilon$ (such may not exist)\footnote{see
Remark 1 to Lemma 5 in \cite{GST08}}.
\label{lem2}
\end{lm}

Different versions of this lemma, as well as similar results
for the flows near a saddle equilibrium state, can be found in
\cite{OSh,GS90,GS92,book}.

In the coordinates of Lemma~\ref{lem1}, we have $P^+_i=(x^+_i,u^+_i,0)$ and $P^-_i=(0,0,y^-_i)\;(i=1,2)$. The neighborhoods $\Pi^+_i$ can be taken as
\begin{equation}\label{eq:glodom}
\begin{array}{l}
\Pi^+_i = \{|x-x^+_i|<\delta^+,\|u-u^+_i\|<\delta^+,|y|<\delta^+\},\\
\Pi^-_i = \{|x|<\delta^-,\|u\|<\delta^-,|y-y^-_i|<\delta^-\}.
\end{array}
\end{equation}
We may also find a convenient form for each global map, using the simplicity of tangencies.

\begin{lm} \label{lem:T_1,2}
\cite{GST08}
If the corresponding homoclinic tangencies are quadratic and simple, then the Taylor expansion (near the point $P_i^-(x=0,u=0,y=y_i^-)$) for map $T_i(\mu_i)\;(i=1,2)$ can be written in  the following form:\\
%in case (1,1) ---
\begin{equation}
\label{eq:T_1,2}
\begin{array}{rcl}
\bar x - x^+_i & = & a_ix + b_i(y-y^-_i) + p_i u +  \dots,\\

\bar u - u^+_i & = & \tilde a_ix + \tilde b_i(y-y^-_i) + \tilde p_i u + \dots,\\
\bar y & = & \mu_i + c_i x + d_i(y-y^-_i)^2 + q_i u +  \dots,
\end{array}
\end{equation}
where $d_i\neq 0$, the dots denote high order terms in the Taylor expansion, and in the saddle case $x\in \mathbb{R}$, $b_i\neq 0,\; c_i \neq 0$ and in the saddle-focus case $x=(x_1,x_2)\in \mathbb{R}^2$, $b_i = (b_{i1},b_{i2})^\top$, $c_i= (c_{i1},c_{i2})$ with $\|b_i\|\neq 0,\; \|c_{i}\| \neq 0$.

All the coefficients shown in
(\ref{eq:T_1,2}) depend on the parameters
$\mu$ (they are at least $C^{r-2}$ in $\mu$). In
particular, $(x^+_i,u^+_i)$ and $y^-_i$ are also
$\varepsilon$-dependent, and at $\mu=0$ they coincide with the
coordinates of the homoclinic points $P^+_i$ and $P^-_i$,
respectively.
\end{lm}

We make some quick comments to this lemma which may be useful for the reader.  First, for $\mu_i=0$, formula (\ref{eq:T_1,2}) shows that $T_i(P_i^-) = P_i^+$. The manifold $W^u_{loc}$ has equation $x=0,y=0$, hence, $T_i{W^u_loc}$ has the equation $\bar x - x^+_i  =   (y-y^-_i) (b_i +  \dots), \bar u - u^+_i  = (y-y^-_i)(\tilde b_i + \dots),
\bar y  =  \mu_i + (y-y^-_i)^2 (d_i +   \dots)$. It follows that the tangency (at $\mu_i=0$) is quadratic if $d_i\neq 0$, and this tangency split generally when varying $\mu_i$.
The equation of the tangent space $\mathcal{T}_{P^-_i}W^{ue}_{loc}$ is $u=0$. At $\mu_i=0$, the plane $DT_i(\mathcal{T}_{P^-_i}W^{ue}_{loc})$ has the equation  $\bar x - x^+_i  =  a_i x + b_i(y-y^-_i), \; \bar u - u^+_i  =  \tilde a_i x + \tilde b_i(y-y^-_i), \; \bar y  =  c_i x $. The condition of simple tangency means that this plane intersects with the leaf $(x=x_i^+, y=0)$ of the foliation $\mathcal{F}^{ss}$ exactly in one point. Thus,
the system $0  =  a_i x + b_i(y-y^-_i), \; 0  =  c_i x $ has $\{x=y=0\}$ as its unique solution. This implies that $\|b_i\|\neq 0, \|c_i\|\neq 0$.

\subsection{Rescaling lemma}

Now we consider the first-return map $T_{km}$ for a $(k,m)$-orbit; that is, for any point $(x_{02},u_{02},y_{02})\in \Pi^+_2$ satisfying $T_{km}(x_{02},u_{02},y_{02})=(\bar x_{02},\bar u_{02},\bar y_{02})$, we have (see Fig. \ref{frtms12})
\begin{equation}\label{eq:T_km-coor-relation}
(x_{02},u_{02},y_{02})
\xmapsto{T^k_0}(x_{11},u_{11},y_{11})
\xmapsto{T_1}(x_{01},u_{01},y_{01})
\xmapsto{T^m_0}(x_{12},u_{12},y_{12})
\xmapsto{T_2}(\bar x_{02},\bar u_{02},\bar y_{02}).
\end{equation}
\begin{figure}[!htb]
\begin{center}
\includegraphics[width=0.6\columnwidth]
{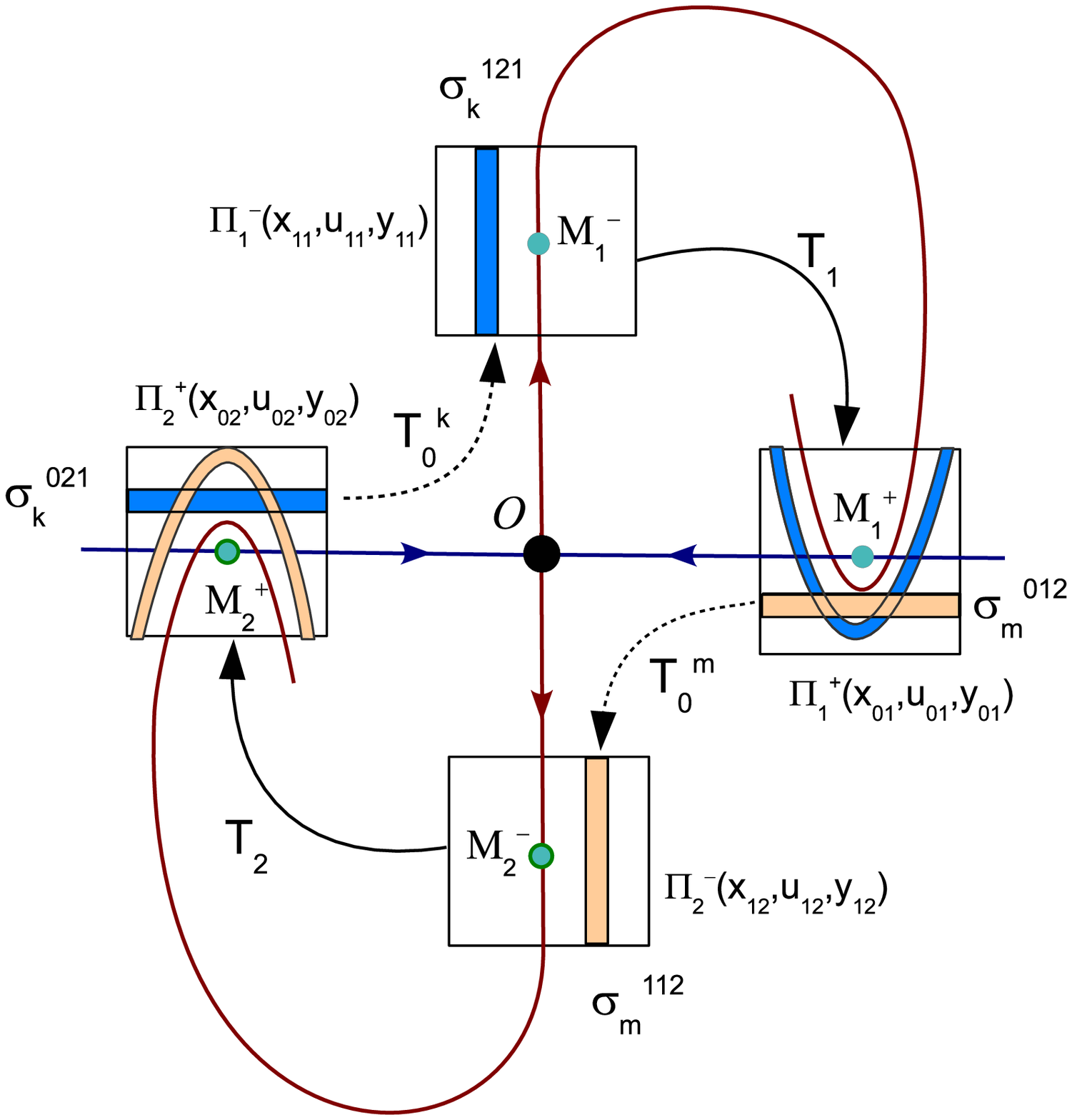} \caption{The action of $T_{km}$ with associated coordinates at some non-zero $\mu$ value.}
\label{frtms12}
\end{center}
\end{figure}
Note that by formula \eqref{T0kk} the ordinate $y_{02}$ is a function of $(x_{02},u_{02},y_{11})$. Hence, we  can study the map $T_{km}$ in the so-called cross-coordinates \cite{Sh67}, i.e. we represent it as $T_{km}:(x_{02},u_{02},y_{11})\mapsto (\bar x_{02},\bar u_{02},\bar y_{11})$.

We first assume that $k\geq m$.
Let $\Pi^+_{i,x,u}$ be the projection of $\Pi^+_i$ to the $x,u$-plane and $\Pi^-_{i,y}$ be the projection of $\Pi^-_i$ to the $y$-axis (see \eqref{eq:glodom}).

\begin{lm} \label{lem:T_km_nf} %\textsf{[Rescaling lemma]}\\
Consider the first return map $T_{km}$ in accordance with relation  \eqref{eq:T_km-coor-relation} with $k\geq m$. There exists a linear coordinate transformation $(x_{02},u_{02},y_{11})\mapsto(X_2,U_2,Y_1)$ such that, in the new coordinates, the map $T_{km}$ takes the form
\begin{equation}\label{eq:T_km_nf:0.1}
\begin{array}{rcl}
\overline X_{21} & =& M_1   - Y_1^2 +  \phi_1(X_2,U_2,Y_1)\;, \\
%  \lambda^m\gamma^k O\left(\gamma^{-(2k+m)/3}\right) + o(1)_{k,m\to+\infty}\;,\\ \\
%
\overline Y_1 & =&
%- \sqrt[3]{d_2d_1^2} \gamma^m \gamma^{\frac{k}{3}}\gamma^{\frac{2m}{3}}( \mu_1 + c_1 \lambda^k\xi_2^+ - \gamma^{-m} y_1^-)
M_2  - \overline  X_{21}^2 + C_2 \lambda^m \gamma^k Y_{1} + \phi_2(X_2,U_2,Y_1)\;, \\
% \lambda^m\gamma^k O\left(\delta_m +  (\hat\lambda/\lambda)^{m} + \gamma^{-(2k+m)/3}\right) +
%o(1)_{k,m\to+\infty} \;, \\ \\
%
(\overline X_{22},U_2) &=  & \phi_3(X_2,U_2,Y_1),
%\lambda^m\gamma^k O\left(\delta_m^{-1}\gamma^{-(2k+m)/3}\right) + O\left(\delta_m^{-1}\gamma^{-(2k+m)/3}\right), \\ \\
%
\end{array}
\end{equation}
where in the saddle case the coordinate ${X}_{22}$ is absent,
$\lambda = \lambda_1$ and $C_2 = b_1c_2$; in the saddle-focus case \\  % $(C_1,C_2)=(b^2_1c_1,b^1_1c_2)$
\begin{equation}\label{eq:T_km_nf:0.1a}
C_{2} = \sqrt{\left(b_{11}^2 + b_{12}^2 \right)\left(c_{21}^2 + c_{22}^2 \right)}\cdot\cos(m\varphi - \nu)
\end{equation}
with $\nu$ given by
$$
\sin\nu = \frac{b_{11}c_{22} - b_{12}c_{21}}{C_2}\quad\mbox{and}\quad  \cos\nu = \frac{b_{11}c_{21} + b_{12}c_{22}}{C_2}.
$$
%
%$C_2 = b_1\sqrt{c_{21}^2 + c_{22}^2} \cos(m\varphi + \beta)$, where $\tan\beta = c_{21}/c_{22}$, in case (2,1);
The rescaled parameters $M_{1,2}$ are given by
\begin{equation}\label{eq:T_km_nf:0.2}
\begin{array}{l}
\displaystyle M_1 = - \sqrt[3]{d_1d_2^2} \gamma^{\frac{4m+2k}{3}} \left( \mu_1 %+ c_1 \lambda^k x_2^+
- \gamma^{-m} y_2^- + O(\lambda^k + \hat\gamma^{-m})\right), \\[5pt]
\displaystyle M_2 = - \sqrt[3]{d_2d_1^2} \gamma^{\frac{4k+2m}{3}} \left( \mu_2 + \alpha_2 \lambda^m  - \gamma^{-k} y_1^-  + O(\hat\lambda^m + \hat\gamma^{-k})\right),
\end{array}
\end{equation}
where in the saddle case $\alpha_1 = c_1 x_2^+$; in the saddle-focus case $\alpha_2 = \tilde C\lambda^m\cos(m\varphi + \tilde\beta)$, where
$$
\tilde C = \sqrt{(c_{21}^2 + c_{22}^2)((x_{11}^+)^2 + (x_{12}^+)^2)}, \;  \sin\tilde\beta = \frac{c_{21}x_{12}^+ - c_{22}x_{11}^+}{\tilde C},\;
\cos\tilde\beta = \frac{c_{21}x_{11}^+ + c_{22}x_{12}^+}{\tilde C}\;.
$$
%, in the case (2,1);
%
The functions $\phi_i$ satisfy
\begin{equation}\label{eq:T_km_nf_nonlinear}
\begin{array}{l}
\phi_1(X_2,U_2,Y_1) =  \lambda^m\gamma^k O\left(\gamma^{-(2k+m)/3}\right) + O(\gamma^{-(2k+m)/3})\;, \\ %o(1)_{k,m\to+\infty}\;, \\
\phi_2(X_2,U_2,Y_1) = \lambda^m\gamma^k O\left(\delta_{km} +  (\hat\lambda/\lambda)^{m} + \gamma^{-(2k+m)/3}\right) +  O(\gamma^{-(2k+m)/3}) + (\hat\gamma/\gamma)^{-k}) \;, \\ %o(1)_{k,m\to+\infty} \;, \\
\phi_3(X_2,U_2,Y_1) = \lambda^m\gamma^k O\left(\delta_{km}^{-1}\gamma^{-(2k+m)/3}\right) + O\left(\delta_{km}^{-1}\gamma^{-(2k+m)/3}\right),
\end{array}
\end{equation}
and their derivatives up to $r$ with respect to variables and up to $(r-2)$ with respect to parameters also satisfy the above estimation. The coefficient $\delta_{km}$ is a freely chosen function of $m$ and $k$ such that $\delta_{km}\to 0$ as $k,m\to +\infty$.

Moreover, when $(x_{02},u_{02},y_{11})$ runs over $\Pi^+_{2,x,u}\times \Pi^-_{1,y}$, the new coordinates cover a centered at zero ball in $\mathbb{R}^{n+1}$ whose radius grows to infinity as $k,m\to \infty$.
\end{lm}
We give the proof of Lemma~\ref{lem:T_km_nf} in
Section~\ref{sec:prlem4}.
Analogous results hold in the case $k<m$. This is done by considering symmetrically for
map $T_{km}$ in the form $T_{km} = T_1 T_0^k T_2 T_0^m $, i.e, in accordance with relation
\begin{equation}\label{eq:T_km-coor-relation2}
(x_{01},u_{01},y_{01})
\xmapsto{T^m_0}(x_{12},u_{12},y_{12})
\xmapsto{T_2}(x_{02},u_{02},y_{02})
\xmapsto{T^k_0}(x_{11},u_{11},y_{11})
\xmapsto{T_1}(\bar x_{01},\bar u_{01},\bar y_{01}),
\end{equation}
when we take the starting point in $\Pi_1^+$,
see Fig.~\ref{frtms12}. Then the proof is conducted in the same way as for Lemma~\ref{lem:T_km_nf} with evident exchanges between coordinates and parameters. Therefore, we formulate the following result without the proof.

\begin{lm} \label{lem:T_km_nf2} %\textsf{[Rescaling lemma]}\\
Consider the first return map $T_{km}$ in accordance with relation \eqref{eq:T_km-coor-relation2} with $k < m$.
There exists a linear coordinate transformation $(x_{01},u_{01},y_{12})\mapsto(X_1,U_1,Y_2)$ such that, in the new coordinates, the map $T_{km}$ takes the form
\begin{equation}\label{eq:T_km_nf:0.11}
\begin{array}{rcl}
\overline X_{11}  &= &M_2   - Y_2^2 +   \psi_1(X_2,U_2,Y_1), \\ % \lambda^k\gamma^m O\left(\gamma^{-(2m+k)/3}\right) + O(\gamma^{-(2m+k)/3})\;, \\
\overline Y_2 & =&
M_1  - \overline  X_{11}^2 + C_1 \lambda^k \gamma^m Y_{2} +  \psi_2(X_2,U_2,Y_1), \\
(\overline X_{12},U_1) &= &\psi_3 (X_2,U_2,Y_1), %\lambda^k\gamma^m O\left(\delta_{km}^{-1}\gamma^{-(2m+k)/3}\right) + O\left(\delta_{km}^{-1}\gamma^{-(2m+k)/3}\right),
\end{array}
\end{equation}
where in the saddle case the coordinate ${X}_{12}$ is absent,
$\lambda = \lambda_1$ and $C_1 = b_2c_1$; in the saddle-focus case \\  % $(C_1,C_2)=(b^2_1c_1,b^1_1c_2)$
\begin{equation}\label{eq:C1}
C_{1} = \sqrt{\left(b_{21}^2 + b_{22}^2 \right)\left(c_{11}^2 + c_{12}^2 \right)}\cdot\cos(k\varphi - \tilde\nu)
\end{equation}
with $\hat\nu$ given by
$$
\sin\hat\nu = \frac{b_{21}c_{12} - b_{22}c_{11}}{C_1}\quad\mbox{and}\quad \cos\hat\nu = \frac{b_{21}c_{11} + b_{22}c_{12}}{C_1}.
$$
%
%$C_2 = b_1\sqrt{c_{21}^2 + c_{22}^2} \cos(m\varphi + \beta)$, where $\tan\beta = c_{21}/c_{22}$, in case (2,1);
The rescaled parameters $M_{1,2}$ are given by
\begin{equation}\label{eq:T_km_nf:0.21}
\begin{array}{l}
\displaystyle M_1 = - \sqrt[3]{d_1d_2^2} \gamma^{\frac{4m+2k}{3}} \left( \mu_1 - \gamma^{-m} y_2^- + \alpha_1 \lambda^k  + O(\hat\lambda^k + \hat\gamma^{-m})\right), \\[5pt]
\displaystyle M_2 = - \sqrt[3]{d_2d_1^2} \gamma^{\frac{4k+2m}{3}} \left( \mu_2  - \gamma^{-k} y_1^-  + O(\lambda^m + \hat\gamma^{-k})\right),
\end{array}
\end{equation}
where in the saddle case $\alpha_1 = c_1 x_2^+$; in the saddle-focus case $\alpha_1 = \hat C\lambda^k\cos(k\varphi + \hat\beta)$, where
$$
\hat C = \sqrt{(c_{11}^2 + c_{12}^2)((x_{21}^+)^2 + (x_{22}^+)^2)}, \;  \sin\hat\beta = \frac{c_{11}x_{22}^+ - c_{12}x_{21}^+}{\hat C},\;
\cos\hat\beta = \frac{c_{11}x_{21}^+ + c_{12}x_{22}^+}{\hat C}\;.
$$
%, in the case (2,1);
%
The functions $\psi_i$ satisfy
\begin{equation*}\label{eq:T_km_psi}
\begin{array}{l}
\psi_1(X_2,U_2,Y_1) =  \lambda^k\gamma^m O\left(\gamma^{-(2m+k)/3}\right) + O(\gamma^{-(2m+k)/3}), \\
\psi_2(X_2,U_2,Y_1) =  \lambda^k\gamma^m O\left(\delta_{km} +  (\hat\lambda/\lambda)^{k} + \gamma^{-(2m+k)/3}\right) +  O(\gamma^{-(2m+k)/3}) + (\hat\gamma/\gamma)^{-m})\;, \\
\psi_3 (X_2,U_2,Y_1) = \lambda^k\gamma^m O\left(\delta_{km}^{-1}\gamma^{-(2m+k)/3}\right) + O\left(\delta_{km}^{-1}\gamma^{-(2m+k)/3}\right),
\end{array}
\end{equation*}
whose derivatives satisfy estimates similar to those for \eqref{eq:T_km_nf_nonlinear}.

Moreover, when $(x_{01},u_{01},y_{12})$ runs over $\Pi^+_{1,x,u}\times \Pi^-_{2,y}$, the new coordinates cover a centered at zero ball in $\mathbb{R}^{n+1}$ whose radius grows to infinity as $k,m\to \infty$,
\end{lm}

\section{Proofs of Theorems 1 and 2.} \label{sec:proofofthm}

\subsection{Proof of Theorem~\ref{thm:shr1}}
We prove the theorem for the $k\geq m$ case and the other follows by a similar argument (we can consider the map $T_{mk} = T_0^kT_1T_0^mT_2$ instead).
We do not assume here that the homoclinic tangencies are simple. Nevertheless, we can consider in $U_0$ local $C^r$-coordinates $(x,y)\in R^{n}\times R$ in which the local stable and unstable invariant manifolds $W^s_{loc}$ and $W^u_{loc}$ of the point $O$ are straightened, i.e., the local map $T_0(\mu)$ has the form:
$$
T_0 : \bar x = A(\mu) x + g(x,y,\mu), \;\; \bar y = \gamma y + h(x,y,\mu),
$$
where the matrix $A$ has the eigenvalues $\lambda_1,...,\lambda_n$ such that $1 > |\lambda_1|\geq ...\geq |\lambda_n|>0$;  and $f(0,y,\mu)\equiv 0, g(x,0,\mu)\equiv 0$. Since $y$ is one-dimensional and $|\lambda_1||\gamma|<1$, we can assume \cite{GST08} that additionally the identities $h(0,y,\mu)\equiv 0, h_y(x,0,\mu)\equiv 0$ are fulfilled. Then we can apply Lemma 2, where the coordinate $u$ is omitted and $\|A\| = \lambda = |\lambda_1| + \epsilon_0$ with arbitrary small constant $\epsilon_0\geq 0$.
Then the
global maps  $T_{1}\equiv f^{N_1} :\Pi_1^-\to\Pi_1^+$ and $T_{2}\equiv
f^{N_2} :\Pi_2^-\to\Pi_2^+$ can be written in the form
$$
\begin{array}{l}
\bar x_{0i} - x_i^+ = F_i(x_{1i}, y_{1i} - y_i^-,\mu), \\
\bar y_{0i} = G_i(x_{1i}, y_{1i} - y_i^-,\mu),  \;\; i=1,2,
\end{array}
$$
where $F_i(0,0,0) = 0, G_i(0,0,0) = 0$, which means that  $T_i(P_i^-) = P_i^+$ at $\mu=0$, and $\partial G_i/\partial y(0,0,0) = 0, \; \partial^2 G_i/\partial y^2(0,0,0) = 2d_i\neq 0$, which means that $W^u(O)$ and $W^s(O)$ have at $\mu=0$ the quadratic homoclinic tangencies at the points of the homoclinic orbits $\Gamma_1$ and $\Gamma_2$. Thus,
we can write the function $G_i$ in the form $G_i = \mu_i + C(x_{1i}) + d_i (y_{1i} - y_i^-)^2 + O(\|x_{1i}\||y_{1i} - y_i^-| + |y_{1i} - y_i^-|^3)$

Now we can write the first-return map $T_{km}: \Pi_2^+ \to \Pi_2^+$, by virtue of the composition (\ref{eq:T_km-coor-relation}), see Fig. \ref{frtms12}), assuming that $k\geq m$. :

$$%\begin{equation}
\begin{array}{l}
x_{01} - x_1^+ = F_1(x_{11}, y_{11} - y_1^-,\mu),\; y_{01} = G_1(x_{11}, y_{11} - y_1^-,\mu),\\
\mbox{where} \;\; x_{11} = \lambda^k \varphi_k^1(x_{02},y_{11},\mu), \; y_{01} = \gamma^{-m}y_{11} + \gamma^{-2m}\varphi_k^2(x_{02},y_{11},\mu), \\ \\
\bar x_{02} - x_2^+ = F_2(x_{12}, y_{12} - y_2^-,\mu),\; \bar y_{02} = G_2(x_{12}, y_{12} - y_2^-,\mu),\\
\mbox{where} \;\; x_{12} = \lambda^m \varphi_m^1(x_{01},\bar y_{11},\mu), \; \bar y_{02} = \gamma^{-k}\bar y_{11} + \gamma^{-2k}\varphi_k^2(\bar x_{02},\bar y_{11},\mu),
\end{array}
$$
Shifting the coordinates  $x_i = x_{0i}-x_i^+,\; y_i = y_{1i} - y_{i}^-$ we bring map $T_{km}$ to the form
\begin{equation}
\label{eq:tkm1}
\begin{array}{l}
x_1 = O(\lambda^k) + O(y_1), \\
\gamma^{-m} y_2 = \hat \mu_1 + d_1 y_1^2 + \lambda^k O(\|x\|) + \gamma^{-2k}O(y_1) + O(y_1^3), \\
\bar x_2 = O(\lambda^m) + O(y_2), \\
\gamma^{-k} \bar y_1 = \hat \mu_2 + d_2 y_2^2 + \lambda^m O(\|x\|) + \gamma^{-2m}O(y_2) + O(y_2^3), \\
\end{array}
\end{equation}
Now we rescale the coordinates as follows
$$
y_1 = \beta_1 Y_1,\; y_2 = \beta_2 Y_2,\; X_1 = \rho_k^{-1}\beta_1 X_1,\; X_2 = \rho_m^{-1}\beta_2 X_2
$$
where
$$%\begin{equation}
\beta_1 = - \frac{1}{\sqrt[3]{d_2d_1^2}}\gamma^{-\frac{k}{3}}\gamma^{-\frac{2m}{3}},\quad\beta_2 = - \frac{1}{\sqrt[3]{d_1d_2^2}}\gamma^{-\frac{m}{3}}\gamma^{-\frac{2k}{3}},
%\label{resc12-2}
$$%\end{equation}
and $\rho_i$ is a function of $i$ which tends to $0$ as $i\to+\infty$. After this, map (\ref{eq:tkm1}) takes the form
\begin{equation}
\label{eq:tkm2}
\begin{array}{l}
X_1 = O(\lambda^k\gamma^{(k+2m)/3} ) + O(\rho_k), \\
Y_2 = M_1 - Y_1^2 + O(\lambda^k\gamma^m\rho_k^{-1}) + O(\gamma^{-(k+2m)/3}) + O(\gamma^{-(2k+m)/3}), \\
\bar X_2 = O(\lambda^m\gamma^{(2k+m)/3} ) + O(\rho_m), \\
\bar Y_1 = M_2 - Y_2^2 + O(\lambda^m\gamma^k\rho_m^{-1}) + O(\gamma^{-(k+2m)/3}) + O(\gamma^{-(2k+m)/3}). \\
\end{array}
\end{equation}

The statement on the ranges of the rescaled coordinates and parameters
is immediate from \eqref{eq:tkm2}, so we only need to prove that the
quantity
\begin{equation}\label{eq:thm1:0}
S_{km}=|\lambda|^m|\gamma|^k
\end{equation}
converges to zero as $k,m\to \infty$. Let $\hat\theta=-\ln|\lambda|/\ln|\gamma|$, hence $|\lambda| = |\gamma|^{-\hat\theta}$.
If $k\geq m$, condition \eqref{mkthet} means that $(\hat\theta - \delta_1)^{-1} < m/k <1$, i.e.  $k<m(\hat\theta -\delta_1)$, with some small $\delta_1>0$ . Then
$$
S_{km} = |\gamma|^{k - m\tilde\theta} \geq |\gamma|^{m(\tilde\theta -\delta_1) - m\tilde\theta} = |\gamma|^{-m\delta_1} \to 0  \;\;\mbox{as}\;\; m\to\infty.
$$

\subsection{Proof of Theorem \ref{thm:shr1-3}.} \label{sec:proofth2}
We only prove for the case where $k\geqslant m$ using Lemma \ref{lem:T_km_nf}; when $k<m$, similar argument applies with using Lemma \ref{lem:T_km_nf2} instead of Lemma \ref{lem:T_km_nf}.

Let us start with the saddle case. By Lemma~\ref{lem:T_km_nf}, we need to find a sequence $\{(k_j,m_j)\}$ of positive integers with $k_j\geqslant m$ and $k_j,m_j\to+\infty$ as $j\to+\infty$ such that, when $\mu_3 = \theta-\theta_0$ varies near 0, the coefficients $C_2\lambda^{m_j}\gamma^{k_j}$ of $Y_1$ in \eqref{eq:T_km_nf:0.1} can cover either $[-s_j, -s_j^{-1}]$ or $[s_j^{-1},s_j]$ for some sequence $\{s_j\}$ of positive numbers tending to infinity.

Evidently, since $C_2 = b_1c_2 \neq 0$ is a constant, it is sufficient to check that the function $S_j = |\lambda|^{m_j}|\gamma|^{k_j}$ as $j\to+\infty$ approximate the infinite interval $(0,+\infty)$ when values of $\theta$ runs over an (arbitrary small) neighborhood of $\theta_0$.

We have for $\ln S_j$ that
\begin{equation}\label{eq:thm2:1}
\theta = \frac{k_j}{m_j} - \frac{\ln S_j}{m_j\ln|\gamma|}.
\end{equation}
Take any sequence $\{s_j\}$ of positive numbers with $s_j\to \infty$ as $j\to \infty$. Let $\theta_j^{1}$ and $\theta_j^{2}$ be the corresponding values of $\theta$ at $S_j = s^{-1}_j$ and $S_j=s_j$. Then
$$
\theta_j^{1,2} = \frac{k_j}{m_j} \mp  \frac{\ln s_j}{m_j\ln|\gamma|}.
$$
So, when $\theta$ runs over the interval $[\theta_j^{1},\theta_j^{2}]$, the values of $S_j$ cover $[s_j^{-1}, s_j]$. Let
$$I_j=[\min(\theta_j^{1},\theta_0),\max(\theta_j^{2},\theta_0)].$$
We will find the desired sequence $\{(k_j,m_j)\}$ if $\mbox{diam}\; I_j \to 0$ as $j\to\infty$. This can be easily achieved by taking $m_j,k_j$ such that
$$\dfrac{k_j}{m_j}\to \theta_0 \quad\mbox{and}\quad \frac{\ln s_j}{m_j\ln|\gamma|} \to 0$$
as $j\to \infty.$

The statement on the parameters $\mu_{1,2}$ and $M_{1,2}$ follows immediately from \eqref{eq:T_km_nf:0.2}. The proof for the saddle case will be completed if the functions $\phi_i$ $(i=1,2,3)$ in \eqref{eq:T_km_nf:0.1} along with their derivatives tend to 0 as $j\to+\infty$. This gives us only some restriction for the exponential growth of $\lambda^m\gamma^k$, which, by \eqref{eq:T_km_nf_nonlinear}, is equivalent to require
$$
\lambda^m\gamma^k  < \max \{ \delta_{km}\gamma^{(2k+m)/3},\; \delta_{km}^{-1},\; (\hat\lambda/\lambda)^{-m} \}.
$$
This inequality can be readily fulfilled by taking appropriate $\delta_{km}$.

We proceed to consider the saddle-focus case. According to \eqref{eq:T_km_nf:0.1a}, the coefficient $C_2\lambda^{m_j}\gamma^{k_j}$ in \eqref{eq:T_km_nf:0.1} now has the form
\begin{equation}\label{eq:thm2:2}
C_j:=C\cos(m_j\varphi - \nu)\lambda^{m_j}\gamma^{k_j},
\end{equation}
where $C=\sqrt{\left(b_{11}^2 + b_{12}^2 \right)\left(c_{21}^2 + c_{22}^2 \right)}$. The idea here is the same as in the saddle case, i.e., we find a sequence $\{(k_j,m_j)\}$ of positive integers with $k_j\geqslant m$ and $k_j,m_j\to+\infty$ as $j\to+\infty$ such that by changing the parameter $\mu_3$ the values of $C_j$ can cover the whole real line in the limit. The difference is that we now vary $\varphi$, instead of $\theta$, near $\varphi_0$ such that the value of $\cos(m_j\varphi - \nu)$ is close to zero, and, at the same time, have $\lambda^{m_j}\gamma^{k_j} \to \infty$ as $k,m\to\infty$.

Take any sequence $\{s_j\}$ of positive numbers with $s_j\to \infty$ as $j\to \infty$. Recall that the range of $\arccos$ is $[0,\pi]$. We solve out $\varphi$ from the equations
$$\cos(m_j\phi^{1,2}_j -\nu)= \pm\dfrac{s_j}{C\lambda^{m_j}\gamma^{k_j}}$$
as
\begin{equation}\label{eq:thm2:3}
\begin{array}{rcl}
\varphi^1_j &=& \dfrac{1}{m_j}\left(\arccos \dfrac{s_j}{C\lambda^{m_j}\gamma^{k_j}}-\nu-\dfrac{2\pi n_j}{m_j}\right),\\ \\
\varphi^2_j &=& \dfrac{1}{m_j}\left(\pi - \arccos \dfrac{s_j}{C\lambda^{m_j}\gamma^{k_j}}-\nu-\dfrac{2\pi n_j}{m_j}\right),
\end{array}
\end{equation}
where $n_j$ can be any integers. Hence, when $\varphi$ runs over
$$I_j=[\min(\varphi_j^{1},\varphi_0),\max(\varphi_j^{2},\varphi_0)],$$
the values of $C$ cover the interval $[-s_j,s_j]$.

We now choose $n_j$ and $m_j$ such that
$$\dfrac{n_j}{m_j}\to \dfrac{\varphi_0}{2\pi}  \quad\mbox{as}\quad j\to\infty,$$
and, then, choose $k_j\geqslant m_j$ accordingly so that
$$\dfrac{s_j}{C\lambda^{m_j}\gamma^{k_j}}\to 0 \quad\mbox{as}\quad j\to\infty. $$
In this way, we find from \eqref{eq:thm2:3} that $\mbox{diam}\; I_j\to 0$ as $j\to\infty$.

Arguing as in the saddle case, one achieves the statement on other parameters and shows that the functions $\phi_i$ $(i=1,2,3)$ converge to $0$ as $j\to\infty$.\qed

\section{Proof of Lemma \ref{lem:T_km_nf}}  \label{sec:prlem4}

As mentioned before, we study the cross-form of the map ${T}_{km}$ according to \eqref{eq:T_km-coor-relation}. Using formulas (\ref{T0kk}) and (\ref{eq:T_1,2})
we can write the relation \eqref{eq:T_km-coor-relation} as
\begin{equation}\label{eq:T_km_nf:1}
\begin{array}{rcl}
(x_{11} - A_1^k x_{02},u_{11}) &=& \hat\lambda^k \left(g_k(x_{02},u_{02},y_{11},\varepsilon),\hat g_k(x_{02},u_{02},y_{11},\varepsilon)\right),\\
\bar y_{02} &=&\gamma^{-k}\bar{y}_{11}+\hat\gamma^{-k}h_k(\bar x_{02},\bar u_{02},\bar{y}_{11},\varepsilon),
\end{array}
\end{equation}
\begin{equation}\label{eq:T_km_nf:2}
\begin{array}{rcl}
x_{01} - x_1^+ &=& a_1 x_{11} + b_1 (y_{11}-y_1^-) +  p_1 u_{11} + \dots,\\
u_{01} - u_1^+ &= &\tilde a_1 x_{11} + \tilde b_1 (y_{11}-y_1^-) + \tilde p_1 u_{11} + \dots,\\
y_{01}& =&  \mu_1 + c_1 x_{11} + d_1 (y_{11}-y_1^-)^2  +  q_1  u_{11} +\dots,
\end{array}
\end{equation}
\begin{equation}\label{eq:T_km_nf:3}
\begin{array}{rcl}
(x_{12} - A_1^m x_{01},u_{12}) &= &\hat\lambda^m \left(g_m(x_{01},u_{01},y_{12},\varepsilon),\hat g_m(x_{01},u_{01},y_{12},\varepsilon)\right),\\
  y_{01}& = &\gamma^{-m} y_{12}  +
\hat\gamma^{-m} h_m(x_{01},u_{01},y_{12},\varepsilon) ,
\end{array}
\end{equation}
\begin{equation}\label{eq:T_km_nf:4}
\begin{array}{rcl}
\bar x_{02} - x_2^+& =& a_2 x_{12} + b_2 (y_{12}-y_2^-) +  p_2 u_{12} + \dots,\\
\bar u_{02} - u_2^+ &=& \tilde a_2 x_{12} + \tilde b_2 (y_{12}-y_2^-) + \tilde p_2 u_{12} +\dots,\\
\bar y_{02}& = & \mu_2 + c_2 x_{12} + d_2 (y_{12}-y_2^-)^2  +  q_2  u_{12} +\dots.
\end{array}
\end{equation}

We construct the sought linear change of coordinates in several steps.

\noindent \textbf{Step 1.} We shift the coordinates with

$$
\begin{array}{l}
x_i = x_{0i} - x_i^+,\; u_i = u_{0i} - u_i^+,\; y_i = y_{1i} - y_i^-,\;\; i=1,2.
\end{array}
$$
Then
$$
\begin{array}{l}
x_{11} =  A_1^k x_{2} +  A_1^k x_{2}^+ + O(\hat\lambda^k), \; u_{11} = O(\hat\lambda^k),\;\; y_{02}= \gamma^{-k} y_1 + \gamma^{-k}y_1^- + O (\tilde\gamma^{-k}), \\
x_{12} =  A_1^m x_{1} +  A_1^m x_{1}^+ + O(\hat\lambda^m), \; u_{12} = O(\hat\lambda^m),\;\; y_{01}= \gamma^{-m} y_2 + \gamma^{-m}y_2^- + O (\tilde\gamma^{-m}),
\end{array}
$$
and the system (\ref{eq:T_km_nf:1})--(\ref{eq:T_km_nf:4}) recasts as (recall that $k\geq m$)
\begin{equation}
\begin{array}{rcl}
x_{1} &=& b_1 y_1 + O(\lambda^k) + O(y_1^2), \\
  u_{1}&=&\tilde b_1 y_1 + O(\lambda^k) + O(y_1^2),\\
\gamma^{-m} y_2 &=& \tilde \mu_1
+ d_1 y_1^2 + \lambda^k O(\|(x_2,u_2)\|) + O\left(|y_1|\|(x_2,u_2)\| + |y_1|^3\right), \\ \\
\bar x_{2} &=& b_2 y_2 + O(\lambda^m) + O(y_2^2), \\
  \bar u_{2} &=& \tilde b_2 y_2 + O(\lambda^m) + O(y_2^2),\\
\gamma^{-k} \bar y_1 &=& \tilde\mu_2
+ d_2 y_2^2 + C_2 A_1^m x_1 +
\hat\lambda^m O(\|(x_1,u_1)\|) + O\left(\lambda^m|y_2|\|(x_1,u_1)\| + |y_2|^3\right), \\ \\
\end{array}
\label{ser1}
\end{equation}
where $C_2 A_1^m x_1 = c_2\lambda^m x_1$ in the saddle case, $C_2 A_1^m x_1 = (c_{21},c_{22}) A_1^m (x_{11},x_{12})^\top$ in the saddle-focus case),
$$
\begin{array}{l}
\tilde \mu_1  = \mu_1 - \gamma^{-m}y_2^- + O(\lambda^k + \hat\gamma^{-m}), \\
\tilde\mu_2 = \mu_2 - \gamma^{-k}y_2^- + \alpha_2\lambda^m + O(\lambda^k + \hat\gamma^{-k} + \hat\lambda^m)
\end{array}
$$
and $\alpha_2$ is the same as in the statement of Lemma~\ref{lem:T_km_nf}.

\noindent \textbf{Step 2.}  The above coordinate transformation gives rise to constant terms in the equations for $x_1,u_1,\bar x_2$ and $\bar u_2$, and to linear-in-$y_1$ and linear-in-$y_2$ terms in equations for $y_2$ and $\bar y_1$, respectively. In order to kill those terms, we shift the coordinates further with
$$
(x_1,u_1,y_1) \mapsto (x_1,u_1,y_1) + O(\lambda^k),\quad (x_2,u_2,y_2) \mapsto (x_2,u_2,y_2) + O(\lambda^m),
$$
where the terms $O(\lambda^k)$ and $O(\lambda^m)$ are chosen appropriately. After that, the formula (\ref{ser1}) can be rewritten as
\begin{equation}
\begin{array}{rcl}
x_{1} &=& b_1 y_1 + \lambda^k O(\|(x_2,u_2)\|) + O(y_1^2), \\
  u_{1} &=& \tilde b_1 y_1 + \lambda^k O(\|(x_2,u_2)\|) + O(y_1^2),\\
\gamma^{-m} \bar y_2 &=&  \tilde\mu_1
+ d_1 y_1^2 + \lambda^k O(\|(x_2,u_2)\|) + O\left(\lambda^k|y_1|\|(x_2,u_2)\| + |y_1|^3\right), \\ \\
\bar x_{2} &=& b_2 y_2 + \lambda^m O(\|(x_1,u_1)\|) + O(y_2^2), \\
\bar u_{2} &=& \tilde b_2 y_2 + \lambda^m O(\|(x_1,u_1)\|) + O(y_2^2),\\
\gamma^{-k} \bar y_1 &=& \tilde\mu_2 +
d_2 y_2^2 + C_2 A_1^m x_1 +
\hat\lambda^m O(\|(x_1,u_1)\|) + O\left(\lambda^m|y_2|\|(x_1,u_1)\| + |y_2|^3\right), \\ \\
\end{array}
\label{ser2}
\end{equation}

\noindent \textbf{Step 3.} We apply affine transformation for $(x,u)$-coordinates.
In the saddle case we put, since $b_1 b_2\neq 0$,
$$
u_i^{new} = u_i - \frac{\tilde b_i}{b_i} x_i, \;\; i=1,2.
$$
In the saddle-focus case, when $x_i = (x_{i1},x_{i2})$ , $b_i = (b_{i1},b_{i2})$ and $b_{i1}^2 + b_{i2}^2 \neq 0$ $(i=1,2)$,  we assume, without loss of generality, that $b_{11}\neq 0$ and $b_{21}\neq 0$.
We put
$$
x_{i2}^{new} = x_{i2} - \frac{b_{i2}}{b_{i1}} x_{i1},\;  u_i^{new} = u_i - \frac{\tilde b_i}{b_{i1}} x_{i1}, \;\; i=1,2.
$$
Then, system (\ref{ser2}) takes the form as
\begin{equation}
\begin{array}{rcl}
x_{11} &=& b_{11} y_1 + \lambda^k O(\|(x_2,u_2)\|) + O(y_1^2),\\
 (x_{12}, u_{1}) &=&  \lambda^k O(\|(x_2,u_2)\|) + O(y_1^2),\\
\gamma^{-m} \bar y_2 &=&  \tilde\mu_1
+ d_1 y_1^2 + \lambda^k O(\|(x_2,u_2)\|) + O\left(\lambda^k|y_1|\|(x_2,u_2)\| + |y_1|^3\right), \\ \\
\bar x_{21} &=& b_{21} y_2 + \lambda^m O(\|(x_1,u_1)\|) + O(y_2^2), \\
(\bar x_{22}, \bar u_{2}) &=&  \lambda^k O(\|(x_1,u_1)\|) + O(y_2^2),\\
\gamma^{-k} \bar y_1 &=& \tilde\mu_2 +
d_2 y_2^2 + \tilde C_2 \lambda^m x_1 +
\hat\lambda^m O(\|(x_1,u_1)\|) + O\left(\lambda^m|y_2|\|(x_1,u_1)\| + |y_2|^3\right), \\ \\
\end{array}
\label{ser3}
\end{equation}
where in the saddle case the coordinates $x_{12}$ and $x_{22}$ are absent, and $x_{11} = x_1, x_{21} = x_2$ and $\tilde C_2 = c_2$;
in the saddle-focus case $\tilde C_2 x_2 = \tilde C_{21} x_{21} + \tilde C_{22} x_{22}$, where
$$
\tilde C_{21} = \frac{1}{b_{21}}\sqrt{\left(b_{21}^2 + b_{22}^2 \right)\left(c_{21}^2 + c_{22}^2 \right)}\cdot\cos(m\varphi - \nu),\;
\nu = \arctan \frac{b_{21}c_{22} - b_{22}c_{21}}{b_{21}c_{21} + b_{22}c_{22}}
$$
and $\tilde C_{22} = \sqrt{c_{21}^2 + c_{22}^2}\cos(m\varphi+\tilde\nu)$, with $\tilde\nu =\arctan (c_{21}/c_{22})$.

\noindent \textbf{{Step 4.} } Now we rescale the coordinates in system (\ref{ser3}) as follows:
\begin{equation}\label{eq:T_km_nf:coor2}
\begin{array}{l}
y_1 = \beta_1 Y_1,\quad x_{11} = b_{11} \beta_1 X_1,\quad (x_{12},u_1) = \tilde\delta_{km} \beta_1 (X_{12},U_1)  \\
y_2 = \beta_2 Y_2,\quad x_{21} = b_{21} \beta_2 X_2,\quad (x_{22},u_2) = \delta_{km} \beta_2 (X_{22},U_2)
\end{array}
\end{equation}
where
\begin{equation}
\beta_1 = - \frac{1}{\sqrt[3]{d_2d_1^2}}\gamma^{-\frac{k}{3}}\gamma^{-\frac{2m}{3}},\quad\beta_2 = - \frac{1}{\sqrt[3]{d_1d_2^2}}\gamma^{-\frac{m}{3}}\gamma^{-\frac{2k}{3}},
\label{resc12-2}
\end{equation}
and $\tilde\delta_{km}$ and $\delta_{km}$ are some functions of $m$ and $k$ which tend to $0$ as $k,m\to+\infty$.

Then system (\ref{ser3}) recasts as
\begin{equation}\label{eq:T_km_nf:6}
\begin{array}{rcl}
X_{11}  &=&  Y_1 +   O\left(\lambda^k\gamma^{-(k-m)/3}\right) +  O\left(\gamma^{-(k+2m)/3}\right), \\[5pt]
(X_{12},U_1) &=&  \tilde\delta_{km}^{-1}O\left(\lambda^k\gamma^{-(k-m)/3}\right) + \delta_{km}^{-1} O\left(\gamma^{-(k+2m)/3}\right), \\[5pt]
Y_2  &=&

M_1   - Y_1^2  +  O\left(\lambda^k\gamma^m\gamma^{-(k-m)/3} + \hat\gamma^{-m}\gamma^m + \gamma^{-(k+2m)/3}\right) , \\ \\

\overline X_{21}  &=&   Y_2 +   O\left(\lambda^m\gamma^k\gamma^{-(2k+m)/3} + \gamma^{-(2k+m)/3} \right),\\ [5pt]
(\overline X_{22},U_2) &=&  \delta_{km}^{-1}O\left(\lambda^m\gamma^k \gamma^{-(2k+m)/3}\right) + \delta_{km}^{-1} O\left(\gamma^{-(2k+m)/3}\right), \\[5pt]
\overline Y_1  &=&
M_2  - Y_2^2  + C_2 \lambda^m \gamma^k X_{11} + \\
&& \qquad + \delta_{km} \lambda^m\gamma^k O(X_{12}) +
 O\left(\hat\lambda^{m}\gamma^k + \hat\gamma^{-k}\gamma^k + \lambda^m\gamma^k\gamma^{-(2k+m)/3}+\gamma^{-(2k+m)/3}\right) ,\\
\end{array}
\end{equation}
where $X_{11},X_{21}\in\mathbb{R}$ and the coordinates $X_{12}$ and $X_{22}$ are absent in the saddle case. Since $k\geq m$, we can take, for example, $\tilde\delta_{km} = \gamma^{-k/3}$.

Finally, expressing intermediate coordinates $(X_1,U_1,Y_2)$ from the first three equations of \eqref{eq:T_km_nf:6} and putting the result into the last three equations, we obtain the following form of $T_{km}$

\begin{equation}\label{ser5}
\begin{array}{rcl}
\overline X_{21}  &=& M_1   - Y_1^2 +   \lambda^m\gamma^k O\left(\gamma^{-(2k+m)/3}\right) + o(1)_{k,m\to+\infty}\;,\\ \\
\overline Y_1  &=&
M_2  - \overline X_{21}^2 + C_2 \lambda^m \gamma^k Y_{1} +  \lambda^m\gamma^k O\left(\delta_{km} +  (\hat\lambda/\lambda)^{m} + \gamma^{-(2k+m)/3}\right) +o(1)_{k,m\to+\infty} \;, \\ \\
(\overline X_{22},U_2) &=& \lambda^m\gamma^k O\left(\delta_{km}^{-1}\gamma^{-(2k+m)/3}\right) + O\left(\delta_{km}^{-1}\gamma^{-(2k+m)/3}\right),
\end{array}
\end{equation}

The last statement of the lemma on the range of the new coordinates follows immediately from \eqref{eq:T_km_nf:coor2} and \eqref{resc12-2}.
\qed

\section{Acknowledgment}

This paper was carried out in the framework of the RSciF
Grant No. 19-11-00280 and partially supported by the  grant 0729-2020-0036 of the Ministry of Science and Higher Education of the Russian Federation (Section 3).  S. Gonchenko thanks the Theoretical Physics and Mathematics Advancement Foundation BASIS, Grant No. 20-7-1-36-5, for support of scientific investigations of his group. D. Li was supported by ERC project 677793 StableChaoticPlanetM.


\begin{thebibliography}{90}
%
\bibitem{Gal93}
J.A.C. Gallas,  Structure of the parameter space of the H\' enon map // Phys. Rev. Lett.  1993, v.70(18), 2714-2717.
%
\bibitem{Gal94}
J.A.C. Gallas, Dissecting Shrimps: Results for Some One-Dimensional Physical Models // Physica A, 1994, 202, 196-223.
%
\bibitem{HGGYK99}
B.R. Hunt, J.A.C. Gallas, C. Grebogi, J.A. Yorke, H. Kocak, Bifurcation rigidity // Physica D, 1999, 129, 35.
%
\bibitem{BGU08}
C. Bonatto, J.A.C. Gallas, Y. Ueda, Chaotic phase similarities and recurrences in a damped-driven D\" uffing oscillator//  Phys. Rev. E, 2008, 77,  026217.
%
\bibitem{Lor08}
E.N. Lorenz, Compound windows of the H\' enon-map // Physica D, 2008, v.237, 1689-1704.
%
\bibitem{G85}
S.V. Gonchenko, On a two parameter family of systems close to a system with a nontransversal
Poincar\'e homoclinic curve. I. // Methods of Qualitative Theory of Differential Equations; Ed.
Gorky St. Univ., 1985, 55-72. English transl. in: Selecta Math. Sovietica, 1990, 10
%
\bibitem{DM00}
H.R. Dullin, J.D. Meiss,  Generalized H\' enon maps: the cubic diffeomorphisms of the plane// Phys. D, 2000, v.143, 262-289.
%
\bibitem{GGO17}
M. Gonchenko, S.V. Gonchenko, I. Ovsyannikov,
Bifurcations of Cubic Homoclinic Tangencies in Two-dimensional
Symplectic Maps // Math. Modeling of Natural Phenomena // 2017, v.12(1).
%
\bibitem{M87}
C. Mir\' a, Chaotic Dynamics. From the One-Dimensional Endomorphism to the Two-dimensional Diffeomorphism // World Scientific, Singapore,1987.
%
\bibitem{BC91a}
J.P. Carcasses, C. Mir\' a, M. Bosh, C. Simo, J.C. Tatjer, ``Crossroad area -- spring area'' transition (I) Parameter plane representation // Bifurcation and Chaos, 1991, v.1, No.1, 183-196.
%
\bibitem{BC91b}
J.P. Carcasses, C. Mir\' a, M. Bosch, C. Simo and J.C. Tatjer. ``Crossroad area -- spring area''
transition (II). Foliated parametric representation // Int. J. Bifur. and Chaos, 1991, 1(2), 339-
348.
%
\bibitem{BC91c}
C. Mir\' a, J.P. Carcasses, ``Crossroad area -- saddle area'' and ``crossroad area -- spring area'' transitions // Bifurcation and Chaos, 1991, v.1, No.3, 641-655.
%
\bibitem{GOT12}
S.V. Gonchenko, I.I. Ovsyannikov, D.V. Turaev, On the effect of invisibility of stable periodic orbits at homoclinic bifurcations // Physica D421, 2012, 1115-1122.
%
\bibitem{book}
L.P. Shilnikov, A.L. Shilnikov, D.V. Turaev, L.O. Chua,
Methods of qualitative theory in nonlinear dynamics //
World Scientific, Part I., 1998; Part II, 2001.
%
\bibitem{N79}
S.E. Newhouse, The abundance of wild hyperbolic sets and non-smooth stable sets
for diffeomorphisms // Publ. Math. Inst. Hautes Etudes Sci., 1979, v.50, 101-151.
%
\bibitem{GST93b}
S.V. Gonchenko, D.V. Turaev, and L.P. Shilnikov, On the existence of Newhouse domains in a neighborhood
of systems with a structurally unstable Poincar\' e homoclinic curve (the higher-dimensional case) // Dokl. Math., 1993,
47 (2), 268-273.
%
\bibitem{PV94}
J. Palis and M. Viana, High dimension diffeomorphisms displaying infinitely many periodic attractors // Ann.
Math., Ser. 2, 1994, 140 (1), 207-250.
%
\bibitem{R95}
N. Romero, Persistence of homoclinic tangencies in higher dimensions // Ergodic Theory Dyn. Syst., 1995, 15 (4),
735-757.
%
\bibitem{GST91a}
S.V. Gonchenko, D.V. Turaev and L.P. Shilnikov,
On models with a structurally unstable homoclinic Poincar\'e
curve // Sov. Math., Dokl. 44, 422-426 (1992).
%
\bibitem{LGYK93}
Lai Y-C., Grebogi C., Yorke J.A., Kan I., How often are chaotic saddles nonhyperbolic // Nonlinearity, 1993, Vol. 6, No. 5, 779-797.
%
\bibitem{GST93a}
S.V. Gonchenko, L.P. Shilnikov and D.V.Turaev,
On models with non-rough Poincar\'e homoclinic curves // Physica D
62, 1-14 (1993).
%
%
\bibitem{Palis2000}
J. Palis, A Global View of Dynamics and a Conjecture on the Denseness of Finitude of Attractors // Asterisque, 2000, v. 261, pp. 339-351.
%
\bibitem{GaS73}
N.K. Gavrilov and L.P. Shilnikov, On  three-dimensional dynamical
systems close to systems with a structurally unstable homoclinic
curve''.- Part~1 // Math.USSR Sb., 1972, v.17, 467-485; Part~2 //
Math.USSR Sb, 1973, v.19, 139-156.
%
\bibitem{N74}
S. E. Newhouse, Diffeomorphisms with infinitely many sinks // Topology, 1974, 13, 9-18.
%
\bibitem{G83}
S.V. Gonchenko, On stable periodic motions in systems that are close to
systems with a structurally unstable homoclinic curve // Rus.
Math. Notes, 1980, v.33, 384-389.
%
\bibitem{GST07}
S.V. Gonchenko, L.P. Shilnikov and D.V.Turaev, Homoclinic tangencies of arbitrarily high orders
in conservative and dissipative two-dimensional maps, {\em Nonlinearity},  2007, {\bf 20}, 241-275.
%
\bibitem{GST93c}
S.V. Gonchenko, D.V. Turaev and L.P. Shilnikov, Dynamical phenomena
in multi-dimensional systems with a non-rough Poincare homoclinic
curve // Russian Acad. Sci.Dokl.Math., 1993, v.47, 3, 410-415.
%
\bibitem{GST08}
S.V. Gonchenko, L.P. Shilnikov and D.V. Turaev, On dynamical
properties of multidimensional diffeomorphisms from Newhouse
regions // \textit{Nonlinearity}, 2008, vol. 21, pp. 923--972.
%
\bibitem{LTY}
L. Tedeshini-Lalli and J.A. Yorke, How often do simple dynamical
processes have infinitely many coexisting sinks // Comm.Math.Phys., 1995, {\bf 106}, pp.635-657.
%
\bibitem{GST96}
S.V. Gonchenko, L.P. Shilnikov and D.V. Turaev, Dynamical phenomena in
systems with structurally unstable Poincar\' e homoclinic orbits // Chaos, 1996,
v. 6, No.1, 15-31.
%
\bibitem{G96} S.V. Gonchenko, Moduli of $\Omega$-conjugacy of two-dimensional
diffeomorphisms with a structurally unstable heteroclinic contour // Sb. Math. 187, 1261-1281 (1996).
%
\bibitem{GST99}
S.V. Gonchenko, D.V. Turaev and L.P. Shilnikov, Homoclinic
tangencies of an arbitrary order in Newhouse domains // Itogi
Nauki Tekh., Ser. Sovrem. Mat. Prilozh. 67, 69-128 (1999) [English
translation in J. Math. Sci. 105, 1738-1778 (2001)].
%
\bibitem{GSV13}
S.V. Gonchenko, C. Simo and A. Vieiro, Richness of dynamics and global bifurcations in systems with a homoclinic figure-eight //
Nonlinearity, 2013, 26(3), 621-678.
%
\bibitem{Be18}
P. Berger, Zoology in the Henon family: twin babies and Milnor's swallows // arXiv:1801.05628v1 [math.DS] 2018.
%
\bibitem{GS87} S.V. Gonchenko and L.P. Shilnikov, Arithmetic properties of topological invariants of systems with structurally unstable homoclinic trajectories // Ukr. Math. J. 39, 15-21 (1987).
%
\bibitem{G89} S.V. Gonchenko, Moduli of systems with non-rough homoclinic
orbits (the cases of diffeomorphisms and vector fields) //
Methods of the qualitative theory and bifurcation theory, 34-49
(1989) [English translation in Sel. Math. Sov. 11, 393-404
(1992)].
%
\bibitem{GS92}
S.V. Gonchenko and L.P. Shil'nikov, On moduli of systems with a
structurally unstable homoclinic Poincar\' e homoclinic curve //
Russian Acad. Sci. Izv. Math., 1993, {\bf 41}, No.3, pp. 417-445.
%
\bibitem{GS90} S.V. Gonchenko and L.P. Shilnikov,
Invariants of $\Omega$-conjugacy of diffeomorphisms with a
nongeneric homoclinic trajectory // Ukr. Math. J. 42, 134-140
(1990).
%
\bibitem{G00}
S.V. Gonchenko, Dynamics and moduli of $\Omega$-conjugacy of
4D-diffeomorphisms with a structurally unstable homoclinic orbit to a
saddle-focus fixed point // Amer. Math. Soc. Transl., 2000, v.200, No.2,
pp.107-134.
%
\bibitem{G02}
S.V. Gonchenko, Homoclinic tangencies, $\Omega$-moduli and
bifurcations // Proc. of the Steklov Inst. of Math., 2002, {\bf
236} .
%
\bibitem{GS86}
S.V. Gonchenko and L.P. Shilnikov, On dynamical systems with
structurally unstable homoclinic curves // Soviet Math. Dokl.,
1986, {\bf 33}, No.1, pp.234-238.
%
\bibitem{GST91b}
S.V. Gonchenko, D.V. Turaev and L.P. Shilnikov,
On models with a non-rough homoclinic Poincar\'e curve // Methods
of qualitative theory and bifurcation theory, Nizhny Novgorod,
36-61 (1991).
%
\bibitem{GStT96}
S.V. Gonchenko, O.V. Sten'kin and D.V. Turaev, Complexity of homoclinic
bifurcations and $\Omega$-moduli // Int.Journal of Bifurcation and Chaos,
v.6, No.6 (1996), pp.969-989.
%
\bibitem{GSStT97}
S.V. Gonchenko, L.P. Shilnikov, O.V. Sten'kin and D.V. Turaev, Bifurcations of systems with structurally unstable homoclinic orbits and $\Omega$-moduli // Computers Math. Applic., 1997,
v. 34, No. 2-4, pp. 111-142.
%
\bibitem{NPT}
S.E. Newhouse, J. Palis and F. Takens, Bifurcations and stability of
families of diffeomorphisms // Publ.Math.Inst. Haute Etudes
Scientifiques, 1983, {\bf 57}, p.5-72.
%
\bibitem{HPS}
M.W. Hirsch, C.C. Pugh and M. Shub, Invariant manifolds // Lecture
Notes in Math., vol.583, Springer-Verlag, Berlin, 1977.
%
\bibitem{T96}
D.V. Turaev, On dimension of non-local bifurcational problems // Int.Journal of Bifurcation and Chaos, 1996, {\bf 6}, No.5,
pp.919-948.
%
\bibitem{OMCMML19}
J. A. de Oliveira, L. T. Montero,  D. R. da Costa,  J. A. Mendez-Bermudez, R. O. Medrano, E. D. Leonel, An investigation of the parameter space for a family of dissipative mappings // Chaos 29, 053114 (2019); https://doi.org/10.1063/1.5048513
%
\bibitem{Afr84}
V.S. Afraimovich, On smooth changes of variables // {\em
Methods of the Qualitative Theory and the Bifurcation Theory},
(E.A.Leontovich-Andronova, ed.), 1984, (Gorky State Univ.),
pp.10-21. (Russian).
%
\bibitem{Sh67}
L.P. Shilnikov, On a Poincar\' e-Birkhoff problem // Math. USSR Sb.,
1967, {\bf 3}, 91-102.
%
\bibitem{Sh68}
L.P. Shilnikov, A contribution to the problem of structure of a
neighborhood of a homoclinic tube of invariant torus // Soviet
Math. Dokl., 1968, {\bf 180}, No.2, 286-289. (Russian).
%
\bibitem{AS73}
V.S. Afraimovich, and L.P. Shilnikov, On critical sets of Morse-Smale
systems // Trans. Moscow Math. Soc., 1973, {\bf 28}.
%
\bibitem{OSh}
I.M. Ovsyannikov and L.P. Shilnikov, On systems with a saddle-focus
homoclinic curve // Math. USSR Sbornik,1987, {\bf 58}, 91-102.
%
\end{thebibliography}
\end{document}